\documentclass[12pt]{article}

\usepackage{amssymb}
\usepackage{graphics} 
\usepackage{url}
\usepackage{hyperref}
\usepackage[authoryear,round]{natbib}
\usepackage{amssymb,color}
\usepackage{graphics,amsmath}

\begin{document}

\newenvironment {proof}{{\noindent\bf Proof.}}{\hfill $\Box$ \medskip}

\newtheorem{theorem}{Theorem}[section]
\newtheorem{lemma}[theorem]{Lemma}
\newtheorem{condition}[theorem]{Condition}
\newtheorem{proposition}[theorem]{Proposition}
\newtheorem{remark}[theorem]{Remark}
\newtheorem{hypothesis}[theorem]{Hypothesis}
\newtheorem{conjecture}[theorem]{Conjecture}
\newtheorem{corollary}[theorem]{Corollary}
\newtheorem{example}[theorem]{Example}
\newtheorem{definition}[theorem]{Definition}

\renewcommand {\theequation}{\arabic{section}.\arabic{equation}}

\def \non {{\nonumber}}
\def \noin {{\noindent}}
\def \hat {\widehat}
\def \tilde {\widetilde}
\def \E {\mathbb{E}}
\def \P {\mathbb{P}}
\def \R {\mathbb{R}}
\def \N {\mathbb{N}}
\def \Z {\mathbb{Z}}
\def \A {\mathbb{A}}
\def \B {\mathbb{B}}

\definecolor{dred}{rgb}{0.83,0.05,0.12}

\title{\large {\bf Obliquely reflecting Brownian motion in nonpolyhedral,}\\
{\bf piecewise smooth cones, with an example of application to}\\ 
{\bf  diffusion approximation of bandwidth sharing queues}}
                                                       
\author{Cristina Costantini\\
Dipartimento di Economia and Unit\`a locale INdAM\\
Universit\`a di Chieti-Pescara\\                           
c.costantini@unich.it}

\date{November 13, 2024}

\maketitle

\begin{abstract}
This work gives sufficient conditions for uniqueness in 
distribution of semimartingale, obliquely reflecting Brownian motion in a 
nonpolyhedral, piecewise ${\cal C}^2$ cone, with radially constant, 
Lipschitz continuous direction of reflection on each face.  

The conditions are shown to be verified by the conjectured 
diffusion approximation to the workload in a class of 
bandwidth sharing networks, thus ensuring that the 
conjectured limit is uniquely characterized. 
This is an essential step in proving the diffusion approximation. 

This uniqueness result is made possible by replacing the 
Krein-Rutman theorem used by Kwon and Williams (1993) 
in a smooth cone with the recent reverse ergodic theorem for 
inhomogeneous, killed Markov chains of Costantini and 
Kurtz (2024). 
\end{abstract}
\vspace{.1in}



\setcounter{equation}{0}

\section{Introduction}\label{sectionintro}
Bandwidth sharing networks operating under a weighted $\alpha$-fair 
bandwidth sharing policy are used to model Internet 
congestion control and were 
studied in depth in \cite{KW04} and \cite{KKLW09}. 
Under this policy, the proportion of the capacity of each 
network resource allocated to each route  
is determined by maximizing a reward function: the 
parameter $\alpha$ is essentially the exponent in the reward function (see 
Section \ref{sectionex} for details). 
The diffusion approximation of the workload process, in the heavy 
traffic regime, is conjectured to 
be a semimartingale obliquely reflecting Brownian 
motion (ORBM) in a piecewise smooth cone. If there are 
only two resources, the cone is a wedge in $\R^2$. If there 
are more than two resources, the 
parameter $\alpha$ determines the shape of the cone: for $\alpha 
=1$, 
the cone is polyhedral; for $\alpha\neq 1$, the cone has curved 
faces that intersect nonsmoothly and can even meet in 
cusp-like singularities (see the example of Section 
\ref{subsecex}).
In the $\alpha =1$ case, \cite{KKLW09} proves the diffusion 
approximation by applying the invariance principle of \cite{KW07} 
and by the uniqueness result of \cite{DW95}. Roughly speaking, 
the invariance principle of \cite{KW07} ensures that the sequence 
of the scaled workload processes is relatively compact and that every 
limit point (in distribution) is a solution of a 
stochastic differential equation with reflection in the 
cone (therefore, as a byproduct, it ensures 
also existence of a solution). Convergence then follows if the solution to the stochastic differential 
equation with reflection is unique in distribution, which, 
for a polyhedral cone, is proved in \cite{DW95}. 
For $\alpha\neq 1$ the authors of \cite{KKLW09} say that their 
argument for the $\alpha =1$ case does not carry over, in particular due to the lack 
of a uniqueness theory for semimartingale ORBM in 
nonpolyhedral cones (see the discussion of the 
literature below) and the conjecture is still open. 

Another notable class of networks in the study of which 
piecewise smooth cones show up is input-queued 
switches under a maximum-weight-$\alpha$ policy. 
Input-queued switches are widely used in Internet 
routing and were studied under a 
maximum-weight-$\alpha$ 
policy in \cite{SW12} and \cite{KW12}. In this 
case the parameter $\alpha$ is the exponent in the criterion 
for determining an optimal matching. As for 
bandwidth sharing networks, the diffusion approximation 
to the workload process is conjectured to 
be a semimartingale obliquely reflecting Brownian 
motion (ORBM) in a piecewise smooth cone, 
which is polyhedral for $\alpha =1$, and has curved 
faces for $\alpha\neq 1$. \cite{KW12} proves the 
diffusion approximation for the $\alpha =1$ case, by 
\cite{KW07} and \cite{DW95}, but for $\alpha\neq 1$ the 
conjecture is still open for the same reasons as for 
bandwidth sharing networks. 

The contribution of this work is precisely to prove, in the nonpolyhedral case, 
a parallel result to \cite{DW95}, that is 
to give sufficient conditions under which uniqueness in 
distribution holds 
for semimartingale ORBM in a nonpolyhedral, piecewise 
smooth cone. The parallelism is incomplete because the 
conditions of \cite{DW95} are also necessary, at least for 
simple polyhedrons, while the conditions given here are 
only sufficient.  
As an example, it will be shown that the conditions 
given here are verified by 
the conjectured diffusion approximation to the workload 
in a class of bandwidth sharing networks. In order to 
complete the proof of the diffusion approximation, one 
needs to show that the  
invariance principle of \cite{KW07} can be applied, or 
prove a suitable result analogous to the invariance 
principle of \cite{KW07}, but this is left for future work. 
These results would also guarantee 
existence of a semimartingale ORBM, so the issue of 
existence is also not treated here. (However see Remark 
\ref{re:exexist} for the existence of a semimartingale 
reflecting ORBM in the example of Section 
\ref{subsecex}.)

Only partial uniqueness results are 
available for obliquely reflecting diffusions in nonsmooth 
domains. 

Most results concern piecewise smooth domains. There are of course the seminal works of 
\cite{HR81} in the orthant and \cite{VW85} and \cite{Wil85} 
in the wedge. The above mentioned paper \cite{DW95} 
proves existence and uniqueness in distribution for semimartingale ORBM 
in a convex polyhedron, with 
constant direction of reflection on each face, 
under a condition on the directions of reflection that generalizes the so 
called completely-${\cal S}$ condition used by  
\cite{TW93} in the orthant. This condition is also necessary for simple polyhedrons. 
Strong existence and uniqueness for a class of ORBMs in 
polyhedral domains that do not satisfy the condition of 
\cite{DW95} are obtained by \cite{Ram06}. In this case the 
ORBM is not a semimartingale. 
For a piecewise 
smooth domain with curved faces and no cusp-like points, 
the best available result is \cite{DI93}, which proves 
strong existence and uniqueness of the solution to a 
stochastic differential equation with reflection, under 
a restriction on the directions of reflection. In dimension 
two, recently \cite{CK24spa} has obtained existence and 
uniqueness in distribution, in a piecewise smooth domain allowing 
for cusps, under a condition that is  
optimal in the sense that it reduces to that of 
\cite{DW95} in the case of a convex polygon with 
constant direction of reflection on each face. This 
condition is also the same under which, in a piecewise 
smooth domain in arbitrary 
dimension, \cite{KR17} prove the equivalence 
between existence and uniqueness in distribution of 
the solution to a stochastic differential equation with 
reflection, well posedness of the corresponding submartingale 
problem and well posedness of the corresponding 
constrained martingale problem (see Section 
\ref{subsecCMP} for constrained martingale problems).  
However \cite{KR17} do not show that the above condition is 
sufficient for existence and uniqueness in distribution. 
In dimension two, the cusp case has been studied by 
\cite{DBT93}, when, on each side of the 
cusp, the direction of reflection forms a constant angle 
with the inward normal, and by \cite{CK18} under the 
same condition as in \cite{CK24spa}. 

However, a piecewise smooth cone is the intersection of 
smooth cones, which are not smooth domains, 
hence it is not a piecewise smooth domain. Therefore 
none of the above results applies to piecewise smooth 
cones. Moreover the arguments used for piecewise 
smooth domains do not extend to piecewise 
smooth cones. 

The work that is most closely related 
to the present one is \cite{KW91}, which proves existence and 
uniqueness of the solution to the submartingale problem 
for an ORBM in a smooth cone, with zero drift and 
radially constant, smooth 
direction of reflection. \cite{KW91} uses smoothness of 
the cone essentially in two points: to solve 
a boundary value problem $(\ref{eq:KW})$), 
the solution of which determines whether the vertex is 
reached; to prove that the transition operator of a 
certain killed Markov chain is a strongly positive 
operator, and hence be able to apply the Krein-Rutman theorem 
(\cite{KR50}) to obtain uniqueness. The Krein-Rutman 
theorem assumes also compactness of the transition 
operator, and proving this requires some delicate oscillation 
estimates. 

The main difference between this work and 
\cite{KW91} is that uniqueness is obtained not by the 
Krein-Rutman theorem, but by a new reverse, ergodic 
theorem for inhomogeneous, killed Markov chains proved in \cite{CK24aap} 
(recalled below as Theorem \ref{th:ergodic}). 
The assumptions of Theorem \ref{th:ergodic} have a clear 
probabilistic meaning and may be verified in a wider 
range of situations: piecewise smooth cones - the case 
of this paper; directions of reflection that are not 
radially constant; variable coefficients of the reflecting 
process; nonsmooth domains that are not cones but can be locally 
approximated by cones (see \cite{CK24aap} in arbitrary 
dimension and \cite{CK24spa} in dimension two), etc.. 
In addition, the assumptions of Theorem 
\ref{th:ergodic} do not involve compactness of the 
transition operators and hence do not 
require oscillation estimates. 

In the present context, Theorem \ref{th:ergodic} allows to 
avoid assuming smoothness of the cone. 
As in \cite{KW91}, here too a function (or two 
functions) are employed to analyze the time till the 
vertex of the cone is reached, but, 
for the purpose of applying Theorem \ref{th:ergodic}, 
the functions do not need to satisfy the 
equalities in $(\ref{eq:KW})$, but only corresponding 
inequalities (Condition \ref{auxfunc}). 
Thus these functions can be found without solving the boundary 
value problem. As an example see the function in Section \ref{subsecex0nothit}. 
Assumption (ii) of Theorem \ref{th:ergodic} is 
related to the strong positivity condition of the 
Krein-Rutman theorem. Here it is verified by a coupling 
result proved in \cite{CK18} and an extension to 
piecewise smooth cones of the support theorem that 
\cite{KW91} uses to prove strong positivity in the smooth 
case (Lemmas \ref{th:scaling} and \ref{th:coupling}). 
This extension is possible thanks to a result 
of \cite{KR17}. A more detailed discussion of how 
the assumptions of Theorem \ref{th:ergodic} are verified 
and a more detailed comparison with \cite{KW91} are 
provided at the beginning of Section 
\ref{sectionpiecesmcone}. 

As in \cite{KW91}, the proof of the main result of this 
paper requires several steps. To 
better illustrate the arguments, first the application to 
bandwidth sharing queues is discussed and a specific 
proof for this example is given: Sections \ref{subsecband} 
and \ref{subsecex} present the model; Section 
\ref{subsecexoutline} provides an outline of the proof; 
the various steps of the proof are then 
carried out each in one of Sections \ref{subsecexDI} to 
\ref{subsecexergodic}. Next the general result is 
presented, discussed and proved in Section 
\ref{sectionpiecesmcone}. 
The assumptions and the main 
result are stated precisely in Section \ref{subsecassum}:  
as in \cite{KW91}, 
it is assumed that the direction of reflection 
on each face is radially constant, but the 
result can be extended to variable directions of reflection (see 
Remark \ref{re:nonradconst}); ${\cal C}^2$ smoothness of each face 
is needed in order to be able to apply the above 
mentioned result by \cite{KR17}. 
The proof is presented following 
the outline of the specific proof for the example: 
the differences between the proof of the general result and that for the 
example are highlighted at the end of Section 
\ref{subsecassum}. The following 
Sections \ref{subseclocaliz} to 
\ref{subsecscale} carry out one step of the proof each. 
Section \ref{subsecCMP} contains some preliminary 
material and the Appendix contains the proof of a more 
technical point.

\subsection{Notation}
For a set $E$, ${\mathbf {1}}_E$ is the indicator function. For a finite set $
E$, $|E|$ denotes its cardinality. 

$\Z_{+}$ denotes the set of nonnegative integers and $\N$ the set 
of strictly positive integers.

$\R_{+}^d$ denotes the nonnegative orthant in $\R^d$ and $\R_{++}^
d$ 
denotes the strictly positive orthant. 

For $e,\,g\in\R^d$, $e\cdot g$ denotes the scalar product. For 
$M\in\R^{d_1\times d_2}$, $M^T$ denotes its transpose $|M|:=\sup_{
x\neq 0}\frac {|Mx|}{|x|}$ 
its norm. 

For $a,\,b\in\R$, $a\wedge b$ denotes the minimum between $a$ and $
b$, while 
$a\vee b$ denotes the maximum. A central dot $\cdot$ sometimes 
replaces the argument of a function. $f\circ g$ denotes the 
composition of two functions. 
For $i,j\in\N$, $\delta_{ij}$ denotes the Kronecker symbol, i.e. $
\delta_{ii}=1$, $\delta_{ij}=0$ for $i\neq j$. 

For $E\subseteq\R^d$, $\overset{\circ}E$ denotes 
the interior of $E$ and $\overline E$ its closure. $d(x,C)$ denotes the 
distance of $x\in\R^d$ from a closed set $C$, i.e. 
$d(x,C):=\inf_{y\in C}|x-y|$. 
${\cal C}^1(C)$ denotes the set of functions continuously 
differentiable on some open neighborhood of $C$, and 
analogously ${\cal C}^2(C)$. ${\cal C}^1_b(C)$ denotes the set of functions in 
${\cal C}^1(C)$ bounded, with bounded derivatives, and analogously 
${\cal C}^2_b(C)$. $S^{d-1}$ denotes the unit sphere in $\R^d$. $
{\cal C}^2(\R^d,\R^d)$ and 
${\cal C}^2(S^{d-1},\R^d)$ denote the set of twice continuously 
differentiable vector fields defined on $\R^d$ and $S^{d-1}$ respectively. 
$B_{\delta}(x)$ denotes the ball of radius $\delta$ centered at $
x$.   
$B_{\delta}$ denotes the ball centered at the origin. 

For $E\subseteq\R^d$, ${\cal C}_{\overline E}[0,\infty )$ denotes the set of 
continuous functions from $[0,\infty )$ to $\overline E$ and ${\cal D}_{\overline 
E}[0,\infty )$ 
denotes the set of functions from $[0,\infty )$ to $\overline E$ that are 
right continuous with left hand limits. For bounded, 
continuous functions $f$, $\|f\|$ denotes the supremum norm. 

For $E\subseteq\R^d$, ${\cal P}(\overline E)$ denotes the set of probability 
measures on $\overline E$. For a random variable $Z$, ${\cal L}(Z
)$ denotes 
the law of $Z$. For a stochastic process $ $$Z$, $\{{\cal F}^Z_t\}$ denotes 
the filtration generated by $Z$, i.e. ${\cal F}^Z_t:=\sigma \{Z(s
),\,s\leq t\}$. 
For a stochastic process $Z$ and a finite random time $\tau$, 
defined on the same probability space, $Z(\tau +\cdot )$ is the 
time shift of $Z$ by $\tau$, that is $Z(\tau +\cdot )$ is defined pathwise 
as 
\[Z(\tau +\cdot )(t):=Z(\tau +t),\quad t\geq 0.\]

\section{An example from diffusion approximation of bandwidth sharing 
queues}\label{sectionex}

\setcounter{equation}{0}

\subsection{Diffusion approximation of bandwidth sharing 
queues}\label{subsecband}

\cite{KW04} and \cite{KKLW09} studied a model of Internet 
congestion control in which the flows present in the 
network share the bandwidth according to a weighted 
$\alpha$-fair policy. 

In their model there is a set $\mathbb{J}=\{1,\cdots ,d\}$ 
of resources and each flow corresponds to the continuous, 
simultaneous transmission of a document through 
a subset of resources that is called a route. The set of 
routes is $\mathbb{I}=\{1,\cdots ,m\}$ and for each $i\in\mathbb{
I}$ route 
$i$ is described by the $i$-th column of an incidence matrix $A$, 
where $A_{ji}:=1$ if route $i$ uses resource $j$ and $A_{ji}:=0$ 
otherwise; it is assumed that for each $j\in\mathbb{J}$ there is at least one 
$i\in\mathbb{I}$ such that $A_{ji}=1$, $A_{li}=0$ for $l\neq j$ 
(Assumption 5.1 in \cite{KKLW09}). 
In particular $A$ has rank $d$. There is no loss of 
generality in assuming, for each $ $$j=1,\cdots ,d$, $A_{ji}=\delta_{
ji}$ for 
$i=1,\cdots ,d$ (while there is no assumption on $A_{j,i}$ for 
$i=d+1,\cdots ,m$, when $m>d$).

At each route $i\in\mathbb{I}$, 
the arrival times of documents are the jump times of a 
Poisson process of parameter $\nu^r_i$, and each document has 
an exponentially distributed size of parameter $\mu^r_i$; 
document sizes are independent of one another and 
independent of all arrival times; the initial number of 
documents is independent of the remaining sizes of the 
documents, which are independent exponentials of 
parameter $\mu^r_i$. All the variables and processes for 
different routes are mutually independent. $r$ is a 
positive, integer scaling 
parameter going to infinity. 

Each resource $j\in\mathbb{J}$ has a capacity $C_j$. 
If $N^r_i(t)$ denotes the number of documents on route $i$ at time 
$t$ and $N^r(t)$ denotes the vector of components $N^r_i(t)$, 
$ $$i\in\mathbb{I}$, for each route $i$ and each resource $j$ 
used by route $i$, the proportion of the capacity $C_j$ allocated to route $
i$ 
at time $t$ is $P_i(N^r(t))$, where the function $P$, which 
represents the allocation policy, is defined in the 
following way. Let $k_i$, $i\in\mathbb{I}$, be strictly 
positive parameters. For each $n\in\R_{+}^m-\{0\}$ let 
$\mathbb{I}_{+}(n):=\{l:n_l>0\}$ and $P^{+}(n)$ be the unique 
maximizer, over all $p\in\R_{+}^{\mathbb{I}_{+}(n)}$ such that 
$\sum_{l\in\mathbb{I}_{+}(n)}A_{jl}p_l\leq C_j$ for all $j\in\mathbb{
J}$, 
of the reward function 
\[F_n(p)=\left\{\begin{array}{ll}
\sum_{l\in\mathbb{I}_{+}(n)}k_ln_l^{\alpha}\,\frac {p_l^{1-\alpha}}{
1-\alpha},&\mbox{\rm if }\alpha\neq 1,\\
\sum_{l\in\mathbb{I}_{+}(n)}k_ln_l\log(p_l),&\mbox{\rm if }\alpha 
=1,\end{array}
\right.\]
where the value of the right member is taken to be $-\infty$ 
if $\alpha\in [1,\infty )$ and $p_l=0$ for some $l\in\mathbb{I}_{
+}(n)$. Then, 
for each $n\in\R_{+}^m-\{0\}$, $P$ is defined as: $P_i(n):=P^{+}_
i(n)$, for 
$i\in\mathbb{I}_{+}(n)$, $P_i(n):=0$ for 
$i\in\mathbb{I}-\mathbb{I}_{+}(n)$; $P(0):=0$. 
An explicit parametric form of the function $P$ is given 
in Proposition 2.1 of \cite{KKLW09} (note that $P$ is 
denoted as $\Lambda$ in \cite{KKLW09}). 

The network defined by the allocation policy $P$ is studied in heavy traffic conditions, that 
is, denoting by $\rho^r_i:=\nu^r_i/\mu^r_i$ and by $\nu^r$, $\mu^
r$ and $\rho^r$ 
the vectors of components $\nu^r_i$, $\mu^r_i$ and $\rho^r_i$, 
$i\in\mathbb{I}$, and by $C$ the 
vector of components $C_j$, $j\in\mathbb{J}$, it is assumed that 
\begin{equation}\nu^r\rightarrow\nu\in\R_{++}^m,\quad\mu^r\rightarrow
\mu\in\R_{++}^m,\quad r\big(A\rho^r-C\big)\rightarrow b\in\R^m,\qquad\mbox{\rm as }
r\rightarrow\infty .\label{eq:heavytr}\end{equation}

Let $M^r:=\mbox{\rm diag}(\mu^r)$, $M:=\mbox{\rm diag}(\mu )$, $\rho_
i:=\nu_i/\mu_i$, $i\in\mathbb{I}$, and 
\begin{equation}\sigma :=2AM^{-1}\mbox{\rm diag}(\nu )M^{-1}A^T.\label{eq:bandsigma}\end{equation}
$\sigma$ is always nonsingular because $A$ has rank $d$. 
$(\ref{eq:heavytr})$ implies 
\begin{equation}A\rho =C.\label{eq:heavytr2}\end{equation}

Building on the fluid model studied in \cite{KW04}, 
\cite{KKLW09} conjecture the following. 

\begin{conjecture}(Conjecture 5.1 of 
\cite{KKLW09})\label{bandconj}
Consider the rescaled workload process 
$X^r$, 
\[X^r(t):=r^{-1}A(M^r)^{-1}N^r(r^2t),\qquad t\geq 0,\]
(denoted as $\hat {W}^r$ in \cite{KKLW09}). Under suitable 
assumptions on $N^r(0)$ and $X^r(0)$, $X^r$ converges in 
distribution, as $r\rightarrow\infty$, to a process $X$ that is a semimartingale Obliquely Reflecting 
Brownian Motion (ORBM), namely a solution of the 
following stochastic differential equation 
with reflection: 
\begin{eqnarray}
&X(t)=X(0)+b\,t+\sigma\,W(t)+\int_0^t\gamma (s)\,d\lambda (s),\quad 
t\geq 0,&\non\\
&X(t)\in\overline {{\cal W}},\quad\lambda (t)=\int_0^t{\bf 1}_{\partial 
{\cal W}}(X(s))d\lambda (s),\quad t\geq 0,&\label{eq:bandSDER}\\
&\quad\gamma (t)\in G(X(t)),\quad |\gamma (t)|=1,\quad d\lambda -
a.e.,\quad t\geq 0,&\non\end{eqnarray}
where ${\cal W}$ is the cone 
\begin{equation}{\cal W}:=\bigg\{x\in\R^d:\,x_j=\sum_{i\in\mathbb{
I}}A_{ji}\frac {\rho_i}{\mu_i(k_i)^{1/\alpha}}\big((q^TA)_i\big)^{
1/\alpha},\;q\in\R^d_{++}\bigg\},\label{eq:alphacone}\end{equation}
\begin{equation}\partial_h{\cal W}:=\{x\in\overline {{\cal W}}:\,
q_h=0,\,q_j>0,\,j\neq h\},\quad h\in\mathbb{J},\label{eq:bandfaces}\end{equation}
are the faces of ${\cal W}$, the direction of reflection on $\partial_
h{\cal W}$ is 
$g^h$: 
\begin{equation}g^h_j:=\delta_{jh},\;j\in\mathbb{J},\label{eq:banddir}\end{equation}
$G(x)$ is the cone of directions of reflection at $x$, 
that is: 
\begin{equation}G(x):=\{g:\,g=\sum_{h:\,x\in\overline {\partial_h
{\cal W}}}u_hg^h,\;u_h\geq 0\},\quad\mbox{\rm for }x\in\partial {\cal W}
-0\label{eq:bandrefcone}\end{equation}
and 
\begin{equation}G(0):=\{g:\,g=\sum_{h\in\mathbb{J}}u_hg^h,\;u_h\geq 
0\}=\R^d_{+}.\label{eq:bandrefcone0}\end{equation}
\end{conjecture}

Note that $\overline {{\cal W}}-\{0\}\subseteq\R_{++}^d$. 

The allocation function $P$ does not enter in the dynamics 
of the limiting process $X$, instead it enters, in its parametric form, in the 
definition of the cone ${\cal W}$: see Theorem 4.1 and $(41)$ 
in \cite{KKLW09}. For $d=2$, the cone ${\cal W}$ is a wedge in 
$\R^2$. For $d>2$, the parameter $\alpha$ plays a crucial role in determining 
the shape of the cone ${\cal W}$: for $\alpha =1$, ${\cal W}$ is a convex 
polyhedral cone; for $\alpha\neq 1$, 
${\cal W}$ is a piecewise smooth cone with curved faces. 
From the example of Section \ref{subsecex}, 
we can get an intuitive idea of the role of $\alpha$: 
for $\alpha >1$, the faces of ${\cal W}$ curve 
outward; for $\alpha <1$, the faces of ${\cal W}$ curve inward and any 
two faces meet in a cusp (Section 5.6 of 
\cite{KKLW09}). 
The parameter $\alpha$ also determines the degree of 
smoothness of each face of the cone. 
In particular the faces are ${\cal C}^2$ smooth if and only if $\alpha\geq 
2$. 

In the $\alpha =1$ case, \cite{KKLW09} proves Conjecture 
\ref{bandconj} by using the results of \cite{KW07} 
(Theorem 5.4, which is an application of the invariance principle of 
Theorem 4.3) and the results of 
\cite{DW95} (Theorem 1.6). Roughly speaking, 
Theorem 4.3 of \cite{KW07} ensures that the sequence 
of processes $\{X^r\}$ is relatively compact and that every 
limit point (in distribution) is a solution of the 
stochastic differential equation with reflection 
$(\ref{eq:bandSDER})$ (therefore, as a byproduct, it ensures 
also existence of a solution). In order 
to have convergence one then needs to 
use Theorem 1.6 of \cite{DW95} which proves uniqueness 
in distribution of the solution to $(\ref{eq:bandSDER})$ when ${\cal W}$ is a polyhedral cone. 

In the $\alpha\neq 1$ case, the authors of \cite{KKLW09} say that 
the above argument cannot be carried out, in particular 
due to the lack of a uniqueness result analogous to that 
of \cite{DW95} for a non-polyhedral cone. In fact, as discussed in the 
Introduction and at the beginning of Section \ref{sectionpiecesmcone}, none 
of the results currently available in the literature applies 
to piecewise smooth,  non-polyhedral cones and there is 
no immediate extension to this case of the results for 
piecewise smooth domains or smooth 
cones (note that a piecewise smooth domain is defined as 
the intersection of a finite number of smooth domains, 
so that, in dimension higher than two, a piecewise smooth cone is not a piecewise smooth 
domain). 

Section \ref{sectionpiecesmcone} of this paper proves, under some 
conditions, that the solution of $(\ref{eq:bandSDER})$ in a 
piecewise ${\cal C}^2$, non-polyhedral cone, is unique in 
distribution, thus providing an essential step towards 
proving the diffusion approximation in the $\alpha\geq 2$ case. 

The example considered in 
Section \ref{subsecex} satisfies all conditions but one for 
arbitrary values of the parameters, and the remaining 
condition for suitable values of the parameters. 
An important feature (which holds more generally for $\alpha >1$) is that, for any 
values of the parameters, any solution of $(\ref{eq:bandSDER})$ 
starting from the origin leaves it immediately with 
probability one. Note however that this is not sufficient to prove uniqueness. 
The arguments used in Section 
\ref{sectionpiecesmcone}, are illustrated below in the 
specific case of the example of Section \ref{subsecex}. 
After describing the example, Section \ref{subsecex} states the 
uniqueness result and the properties of the solution precisely; 
Section \ref{subsecexoutline} 
provides an outline of the proof of the uniqueness result; the main points of the 
proof are then dealt with each in a separate Section.

\subsection{The example}\label{subsecex}

Consider the example where there are $d$ resources, 
each resource has a route that only goes through that 
resource and there is one additional route that goes 
through all resources, that is: $\mathbb{J}:=\{1,\cdots ,d\}$, 
$\mathbb{I}:=\{1,\cdots ,d+1\}$ and 
\begin{equation}A_{ji}=\delta_{ji}\mbox{\rm \ for }j,i=1,\cdots ,
d,\qquad A_{j(d+1)}=1,\mbox{\rm \ for }j=1,\cdots ,d.\label{eq:exincid}\end{equation}
We will take 
\begin{equation}\nu_i=k_i=1,\mbox{\rm \ for }i=,\cdots ,d+1,\qquad
\mu_i=1,\mbox{\rm \ for }i=1,\cdots ,d,\quad\mu_{d+1}=\mu ,\label{eq:exparam}\end{equation}
but the values of these parameters are relevant only in 
Section \ref{subsecex0nothit}. In addition, for simplicity 
we take $d=3$, but everything works in the higher dimensional case. 
Then the cone ${\cal W}$, defined by by $(\ref{eq:alphacone})$, and 
the faces $\partial_h{\cal W}$, defined by $(\ref{eq:bandfaces})$, take the form 
\begin{equation}{\cal W}=\bigg\{x\in\R^3_{++}:\,x_j=q_j^{1/\alpha}
+\frac 1{\mu^2}\big(q_1+q_2+q_3\big)^{1/\alpha},\;q\in\R_{++}^3\bigg
\},\label{eq:excone}\end{equation}
\begin{equation}\partial_h{\cal W}=\bigg\{x\in\R^{3+}_{++}:\,x_h=\frac 
1{\mu^2}\big(\sum_{l\neq h}q_l\big)^{1/\alpha},\,x_j=q_j^{1/\alpha}
+x_h,\;q_j\in\R_{++},\;j\neq h\bigg\}.\label{eq:exfaces}\end{equation}
Straightforward computations show that each face 
$\partial_h{\cal W}$ is ${\cal C}^1$ smooth at each point if and only if $
\alpha >1$, and 
${\cal C}^2$ smooth if and only if $\alpha\geq 2$. We assume $\alpha
\geq 2$,  
but this degree of smoothness is needed only in Section 
\ref{subsecexergodic} to be able to apply Theorem 
\ref{th:ergodic}: in fact, in order to verify the assumptions 
of Theorem \ref{th:ergodic}, one needs to show that, 
starting at a point different from the origin, any 
solution of $(\ref{eq:bandSDER})$ spends zero time on 
$\partial {\cal W}-\{0\}$ before hitting the origin, almost surely. This 
does not follow from the fact that $\partial {\cal W}-\{0\}$ has surface 
measure zero, but is ensured by the ${\cal C}^2$-smoothness of ${\cal W}_
h$ and conditions on the 
directions of reflection that are satisfied in this example 
(see Proposition 2.12 of \cite{KR17} and Section \ref{subsecexDI}). 
The other intermediate results, in particular the fact that 
any solution of $(\ref{eq:bandSDER})$ 
starting from the origin leaves it immediately with 
probability one, hold for $\alpha >1$. 

From $(\ref{eq:exfaces})$ one can compute the unit inward normal 
vector $n^h(x)$ at $x\in\partial_h{\cal W}$: 
\begin{eqnarray}
\label{eq:exnorm}n^h_j(x)&=&-\,c^h_n(q)\frac 1{\mu^2}\,\big(\sum_{
l\neq h}q_l\big)^{(1-\alpha )/\alpha}\,q_j^{(\alpha -1)/\alpha},\qquad 
j\neq h,\non\\
\\
n^h_h(x)&=&\,c^h_n(q)\big[1+\frac 1{\mu^2}\big(\sum_{l\neq h}q_l\big
)^{(1-\alpha )/\alpha}\sum_{l\neq h}\,\,q_l^{(\alpha -1)/\alpha}\big
],\non\end{eqnarray}
where $c^h_n$ is the normalization constant, 
\[c^h_n(q):=\bigg\{\frac 1{\mu^4}\bigg(\sum_{l\neq h}q_l\bigg)^{\frac {
2(1-\alpha )}{\alpha}}\,\sum_{j\neq h}q_j^{\frac {2(\alpha -1)}{\alpha}}\,
+\bigg[1+\frac 1{\mu^2}\bigg(\sum_{l\neq h}q_l\bigg)^{\frac {1-\alpha}{
\alpha}}\sum_{l\neq h}\,\,q_l^{\frac {\alpha -1}{\alpha}}\bigg]^2\bigg
\}^{-\frac 12}\]
Since $\sum_{l\neq h}q_l>0$ for all $x\in\overline {\partial_h{\cal W}}
-\{0\}$, $0<c^h_n(q)\leq 1$ for all 
$x\in\overline {\partial_h{\cal W}}-\{0\}$ and $(\ref{eq:exnorm})$ extends to all $
x\in\overline {\partial_h{\cal W}}-\{0\}$. 

The matrix $\sigma$, defined by $(\ref{eq:exsigma})$, in this case 
takes the form  
\begin{equation}\sigma =2\;\left[\begin{array}{cccc}
\big(1+\frac 1{\mu^2}\big)&\frac 1{\mu^2}&\frac 1{\mu^2}\\
\frac 1{\mu^2}&\big(1+\frac 1{\mu^2}\big)&\frac 1{\mu^2}\\
\frac 1{\mu^2}&\frac 1{\mu^2}&\big(1+\frac 1{\mu^2}\big)\end{array}
\right].\label{eq:exsigma}\end{equation}

For this example, Theorem \ref{th:exuniq} below determines a set of values of 
the parameter $\mu$ for which the solution to 
$(\ref{eq:bandSDER})$ is unique. Here a solution is meant 
in the weak sense and uniqueness is meant in 
distribution. More precisely:

\begin{definition}\label{def:exSDER}
A continuous process $X$ is a solution of 
$(\ref{eq:bandSDER})$ 
if there exist a standard Brownian motion $W$, a continuous, 
nondecreasing process $\lambda$ and a process 
$\gamma$ with measurable paths, all defined on the same probability space as $
X$, 
such that $W(t+\cdot )-W(t)$ is independent of ${\cal F}_t^{X,W,\lambda 
,\gamma}$ for all $t\geq 0$ 
and $(\ref{eq:bandSDER})$ is satisfied almost surely. 
If $X(0)=x$ almost surely for 
some $x\in\overline {{\cal W}}$, we say that $X$ starts at $x$. 

Given an initial distribution on $\overline {{\cal W}}$, {\em uniqueness in distribution\/} holds if any 
two solutions of $(\ref{eq:bandSDER})$ with
the same initial distribution have the same distribution on 
$C_{\overline {{\cal W}}}[0,\infty )$. 
\end{definition}

For a solution $X$ of $(\ref{eq:bandSDER})$, let 
\begin{equation}\vartheta^X:=\inf\{t\geq 0:\,X(t)=0\},\label{eq:ex0hittime}\end{equation}
\begin{equation}\tau^{X,\delta}:=\inf\{t\geq 0:\,|X(t)|=\delta \}
,\quad\delta >0.\label{eq:exdeltahittime}\end{equation}
(When there is no risk of confusion, the superscript $X$ 
will be omitted.)

\begin{theorem}\label{th:exuniq}
If $\mu >\sqrt {\frac 3{\sqrt {2}-1}}$, uniqueness in distribution holds for 
$(\ref{eq:bandSDER})$. 

If $X$ starts at $x\in\overline {{\cal W}}-\{0\}$, 
\begin{equation}\P (\vartheta^X<\infty )=0.\label{eq:ex0nothit}\end{equation}

If $X$ starts at $0$, 
\begin{equation}\lim_{\delta\rightarrow 0}\tau^{X,\delta}=0,\qquad 
a.s.\label{eq:exdeltahit}\end{equation}
\end{theorem}

\begin{remark}\label{re:exlowerbd}
It may seem counterintuitive that 
the solution does not hit the origin 
for large enough values of $\mu$, as large values of $\mu$ 
correspond to a small expected length of the document 
that uses all resources, but one has to take into account 
that, in the heavy traffic regime, $(\ref{eq:heavytr2})$ 
must be satisfied, so that, in some sense, to 
a smaller expected length of the document 
that uses all resources correspond smaller capacities. 
\end{remark}

\begin{remark}\label{re:exexist}
The issue of the existence of solutions to 
$(\ref{eq:bandSDER})$ is not addressed in this paper. 
However it is worth mentioning that for the example of 
this section existence of solutions can be proved by the 
same arguments as in Theorems 3.1 and 4.1 of 
\cite{CK18}. 
\end{remark}
\vskip.3in

\subsection{Outline of the proof of Theorem \ref{th:exuniq}}\label{subsecexoutline}

The proof of Theorem \ref{th:exuniq} requires several 
steps. This section provides an outline of these steps. 
For each step, the key points are proved in one of the following 
Sections. Recall that $B_{\delta}$ denotes the ball of radius 
$\delta$ centered at the origin. 
\vskip.2in

\noindent
{\bf Step 1.} By the localization results of \cite{CK24spa}, it is 
enough to prove uniqueness for $(\ref{eq:bandSDER})$ in each 
of the bounded domains $\big\{{\cal W}_n\big\}_{n\geq 0}$, where $
{\cal W}_0:={\cal W}\cap B_4$, 
${\cal W}_n:={\cal W}\cap\big(\overline {B_{1+2(n-1)}}\big)^c\cap 
B_{6+2(n-1)}$, $n\geq 1$, taking the unit 
inward normal as the direction of reflection on the parts of 
the boundary of ${\cal W}_n$ that are not subsets of the boundary 
of ${\cal W}$. 

Each domain of the form 
${\cal W}\cap\big(\overline {B_{\delta}}\big)^c\cap B_{\delta'}$, $
0<\delta <\delta'$, is piecewise ${\cal C}^2$, 
therefore, taking the unit 
inward normal as the direction of reflection on ${\cal W}\cap\partial 
B_{\delta}$ 
and on ${\cal W}\cap\partial B_{\delta'}$, uniqueness is ensured by 
\cite{DI93} (Section \ref{subsecexDI}). 
In particular uniqueness holds for each ${\cal W}_n$, $n\geq 1$, 
therefore we only need to prove uniqueness in ${\cal W}_0$. 

In addition, the results of \cite{CK19} ensure that there 
exist strong Markov solutions of $(\ref{eq:bandSDER})$ in 
${\cal W}_0$ and that uniqueness among strong Markov solutions is 
equivalent to uniqueness among all solutions. 
Therefore in the sequel only 
strong Markov solutions are considered. 
\vskip.2in

\noindent
{\bf Step 2.} Consider a solution $X$ of $(\ref{eq:bandSDER})$ 
in ${\cal W}_0$ starting off the origin, that is with initial distribution 
$\nu$ with compact support in $\overline {{\cal W}_0}-\{0\}$. 
For all $\delta$ small enough that $B_{\delta}$ has empty intersection 
with the support of $\nu$, the localization results of \cite{CK24spa}, 
imply that $X$ must 
agree with some solution of $(\ref{eq:bandSDER})$ in 
${\cal W}_0\cap\big(\overline {B_{\delta /2}}\big)^c$ 
(with direction of reflection on ${\cal W}_0\cap\partial B_{\delta 
/2}$ the inward 
normal)
up to the first time it hits $\partial B_{\delta}$. Since \cite{DI93} 
guarantees that uniqueness 
holds for solutions of $(\ref{eq:bandSDER})$ in ${\cal W}_0\cap\big
(\overline {B_{\delta /2}}\big)^c$, 
$X$ is uniquely determined up to the first time it hits $\partial 
B_{\delta}$ for 
all $\delta$ sufficiently small, hence it is uniquely determined up to the first time 
it hits the origin. 

Moreover every solution of $(\ref{eq:bandSDER})$ in ${\cal W}_0$ starting 
off the origin with probability one never 
reaches the origin  
(Section \ref{subsecex0nothit}), hence the solution is 
uniquely determined for all times. In particular, for 
every $x\in\overline {{\cal W}_0}-\{0\}$, any two solutions starting at $
x$ have 
the same distribution. 

Note that the fact that, starting off the origin, the 
origin is never hit holds independently of the drift $b$: 
this is because very close to the origin the contribution 
of the drift is negligible  with respect to that of the 
diffusion term and of the reflection (see 
$(\ref{eq:exAV})$). 
\vskip.2in

\noindent
{\bf Step 3.} Let $X$ be an arbitrary solution of $(\ref{eq:bandSDER})$ 
in ${\cal W}_0$ starting at the origin. $X$ immediately leaves the origin, that is 
$(\ref{eq:exdeltahit})$ holds. In fact the stronger 
statement 
\[\lim_{\delta\rightarrow 0}\E[\tau^{X,\delta}]=0.\]
holds (Section \ref{subsecexleave0}). 
\vskip.2in

\noindent
{\bf Step 4.} A key consequence of Step 3 is the following: if 
any two solutions of $(\ref{eq:bandSDER})$ in ${\cal W}_0$ starting 
at the origin, $X$ and $\tilde {X}$, have the same hitting 
distributions on $\partial B_{\delta}$, that is ${\cal L}(X(\tau^{
X,\delta}))={\cal L}(\tilde {X}(\tau^{\tilde {X},\delta}))$, for 
all $\delta >0$ sufficiently small, then any two solutions of $(\ref{eq:bandSDER})$ 
starting at the origin have the same one-dimensional 
distributions. Combined with Step 2, this yields that any 
two solutions with the same initial distribution have the 
same one-dimensional distributions, which gives 
uniqueness of their distribution on the path space 
(Section \ref{subsecexhitdist}). 
\vskip.2in

\noindent
{\bf Step 5.} By the previous step, we are reduced to proving 
that, for all $\delta >0$ sufficiently small, 
for any solution $X$ of $(\ref{eq:bandSDER})$ in ${\cal W}_0$ starting 
at the origin, the hitting 
distribution on $\partial B_{\delta}$ is always the same. 
Let $\tau^{X,2^{-2l}\delta}$ be the first time $X$ hits 
$\partial B_{2^{-2l}\delta}$, $l\in\Z_{+}$. In the sequel of this step, in order to 
make formulas more readable, the following shorthand 
notation: 
\[\tau^l:=\tau^{X,2^{-2l}\delta}\]
is used. 
Note that, since $X$ starts at the origin, 
$\tau^l<\tau^{l-1}$. $X(\tau^l+\cdot )$ is a solution of 
$(\ref{eq:bandSDER})$ with initial distribution supported 
in $\partial B_{2^{-2l}\delta}$, and hence, by Step 2, its distribution is the same for any $
X$.
In particular, since 
\[X(\tau^{l-1})=X(\tau^l+\tau^{X(\tau^l+\cdot ),2^{-2(l-1)}\delta}
),\]
we see that ${\cal L}\big(X(\tau^{l-1})\big|X(\tau^l)=x\big)$, $x
\in\partial B_{2^{-2l}\delta}$, is the same for any $X$. 

Since $X$ is a strong Markov process, for $n\in\N$ 
$\big\{X_k\big\}_{0\leq k\leq n}:=\big\{X(\tau^{n-k})\big\}_{0\leq 
k\leq n}$ is a (inhomogeneous) 
Markov chain and $X_n=X(\tau^{X,\delta})$. Moreover, 
since the origin is never hit by a 
solution of $(\ref{eq:bandSDER})$ between $\tau^{n-k}$ and $\tau^{
n-k-1}$,
by Step 2 two Markov chains $\big\{X_k\big\}_{0\leq k\leq n}$ and $\big
\{\tilde {X}_k\big\}_{0\leq k\leq n}$ corresponding 
to two solutions $X$ and $\tilde {X}$ have the same transition 
kernels and differ only by their initial 
distributions. Then, by the reverse ergodic theorem of 
\cite{CK24aap} (recalled below as Theorem 
\ref{th:ergodic}), ${\cal L}(X_n)$ and ${\cal L}(\tilde {X}_n)$ must converge to the 
same limit as $n\rightarrow\infty$, that is ${\cal L}(X(\tau^{X,\delta}
))={\cal L}(\tilde {X}(\tau^{\tilde {X},\delta}))$ 
(Section \ref{subsecexergodic}). 

\begin{remark}\label{re:exkilled}
As mentioned above, in this example the probability transition kernels 
\[Q_{n-k}(x,B):=\P\big(X(\tau^{n-k-1})\in B\big|X(\tau^{n-k})=x\big
),\]
are the same for all solutions $X$, because 
the origin is never hit by a 
solution of $(\ref{eq:bandSDER})$ between $\tau^{n-k}$ and $\tau^{
n-k-1}$,

If the solutions of $(\ref{eq:bandSDER})$ could hit the origin  
between $\tau^{n-k}$ and $\tau^{n-k-1}$ with positive probability, 
in order to ensure a priori that the transition kernels 
are the same for all solutions, 
one would have to consider the killed 
Markov chains with subprobability transition kernels $ $
\[Q_{n-k}(x,B):=\P\big(\tau^{n-k-1}<\vartheta^{n-k},\,X(\tau^{n-k
-1})\in B\big|X(\tau^{n-k})=x\big),\]
where $\vartheta^{n-k}$ is the first time $X$ hits the origin 
after $\tau^{n-k}$. 

However, the reverse ergodic theorem of \cite{CK24aap} applies to 
killed Markov chains as well, under suitable assumptions 
(see Section \ref{subsecergodic}).
\end{remark}
\vskip.2in

For the proof of $(\ref{eq:ex0nothit})$ and 
$(\ref{eq:exdeltahit})$, see Sections 
\ref{subsecex0nothit} and \ref{subsecexleave0}. 
\vskip.3in

\subsection{Uniqueness in a domain ${\cal W}\cap\big(\overline {B_{\delta}}\big)^
c\cap B_{\delta'}$}\label{subsecexDI}
Let us verify that in each domain of the form 
${\cal W}_{\delta ,\delta'}:={\cal W}\cap\overline B_{\delta}^c\cap 
B_{\delta'}$, $0<\delta <\delta'$, taking the unit 
inward normal as the direction of reflection on $\partial B_{\delta}
\cap {\cal W}$ 
and on $\partial B_{\delta'}\cap {\cal W}$, the assumptions of \cite{DI93} 
are satisfied. These assumptions are satisfied in 
particular if, at each point $x\in\partial {\cal W}_{\delta ,\delta'}$, the set of normal 
vectors is linearly independent and a certain matrix 
(the matrix defined by $(\ref{eq:exDImx})$) has spectral radius strictly less than 
$1$. Denote 
\[\partial_0{\cal W}_{\delta ,\delta'}:=\partial B_{\delta}\cap {\cal W}
,\;\partial_4{\cal W}_{\delta ,\delta'}:=\partial B_{\delta'}\cap 
{\cal W},\;\partial_h{\cal W}_{\delta ,\delta'}:=\partial_h{\cal W}
\cap\big(\overline {B_{\delta}}\big)^c\cap B_{\delta'}\;h\in \{1,
2,3\},\]
\[g^0(x):=n^0(x):=\frac x{|x|},\quad g^4(x):=n^4(x);=-\frac x{|x|}
,\quad g^h(x):=g^h,\;\;h\in \{1,2,3\},\]
\[{\cal J}(x):=\big\{h\in \{0,1,2,3,4\}:\,x\in\overline {\partial_
h{\cal W}_{\delta ,\delta'}}\big\}.\]
Note that, by $(\ref{eq:exnorm})$, 
\[g^h\cdot n^h(x)=n^h_h(x)>0,\quad\mbox{\rm for }x\in\overline {\partial_
h{\cal W}_{\delta ,\delta'}},\;\;h\in \{1,2,3\},\]
and obviously the same holds for $h\in \{0,4\}$. 
Let us check that, at each point on the boundary, all 
normal vectors are linearly independent. 
Clearly, 
\begin{equation}\begin{array}{c}
n^0(x)\cdot n^h(x)=0\qquad\forall x\in\overline {\partial_h\,{\cal W}_{
\delta ,\delta'}}\cap\overline {\partial_0{\cal W}_{\delta ,\delta'}}
,\quad h\in \{1,2,3\},\\
n^4(x)\cdot n^h(x)=0\qquad\forall x\in\overline {\partial_h\,{\cal W}_{
\delta ,\delta'}}\cap\overline {\partial_4{\cal W}_{\delta ,\delta'}}
,\quad h\in \{1,2,3\}.\end{array}
\label{eq:n04}\end{equation}
Let $x\in\overline {\partial_h{\cal W}_{\delta ,\delta'}}\cap\overline {
\partial_k{\cal W}_{\delta ,\delta'}}\cap {\cal W}_{\delta ,\delta'}$, $
h,k\in \{1,2,3\}$, $h\neq k$. Then, by 
$(\ref{eq:exnorm})$, 
\begin{equation}n^h_h(x)>0,\quad n^h_k(x)=0,\label{eq:nh}\end{equation}
so that, of course, $n^h(x)$ and $n^k(x)$ are linearly 
independent. Finally at a point $x\in\overline {\partial_h\,{\cal W}_{
\delta ,\delta'}}\cap\overline {\partial_k{\cal W}_{\delta ,\delta'}}
\cap$$\overline {\partial_0{\cal W}_{\delta ,\delta'}}$, 
$h,k\in \{1,2,3\}$, $h\neq k$, $n^0(x)$ 
is orthogonal to both $n^h(x)$ and $n^k(x)$, which are linearly 
independent among themselves, so the three vectors are 
linearly independent, and similarly at 
$x\in\overline {\partial_h{\cal W}_{\delta ,\delta'}}\cap\overline {
\partial_k{\cal W}_{\delta ,\delta'}}\cap\overline {\partial_4{\cal W}_{
\delta ,\delta'}}$. 

The main assumption of \cite{DI93} is $(3.7)$ of 
\cite{DI93}. A sufficient condition is the condition stated in Remark 3.1   
of \cite{DI93}, namely that, at every point $x\in\partial {\cal W}_{
\delta ,\delta'}$, 
${\cal J}(x)=\{h_1,\cdots h_k\}$, the spectral radius of the matrix of elements 
\begin{equation}\frac {|g^{h_i}(x)\cdot n^{h_j}(x)|}{g^{h_i}(x)\cdot n^{h_i}(x)}
-\delta_{h_ih_j},\quad i,j=1,\cdots ,k,\label{eq:exDImx}\end{equation}
is strictly less than $1$. But, for each $h\in {\cal J}(x)$, 
\[g^0(x)\cdot n^h(x)=n^0(x)\cdot n^h(x),\quad g^4(x)\cdot n^h(x)=
n^4(x)\cdot n^h(x),\quad g^k(x)\cdot n^h(x)=n_k^h(x),\]
and we see, by $(\ref{eq:n04})$ and $(\ref{eq:nh})$, 
that $0$ is the only eigenvalue of the above matrix, 
so the condition is trivially satisfied. 
\vskip.3in

\subsection{Starting off the origin, the origin is never 
reached}\label{subsecex0nothit}
Let $\mathbb{A}$ be the operator 
\begin{equation}\mathbb{A}f(x):=b\cdot\nabla f(x)+\frac 12\mbox{\rm tr}
((\sigma\sigma^T)D^2f(x)),\label{eq:exA}\end{equation}
and consider the function 
\begin{equation}V(x):=|x|^{\beta},\quad x\in\R^3-\{0\},\qquad\mbox{\rm with}
\quad 1-\frac 2{\big(1+\frac 3{\mu^2}\big)^2}<\beta <0.\label{eq:beta}\end{equation}
Recalling $(\ref{eq:bandrefcone})$ and 
$(\ref{eq:bandrefcone0})$ we see that $V$ satisfies 
\[\nabla V(x)\cdot g\leq 0,\quad\forall g\in G(x),\;\;\forall x\in
\partial {\cal W}_0-\{0\},\;|x|<4.\]
In addition, straightforward computations give 
\begin{equation}\mathbb{A}V(x)=\beta |x|^{\beta -2}\bigg\{b\cdot x+\frac 12\mbox{\rm tr}
(\sigma\sigma^T)+\frac 12\frac {\beta -2}{|x|^2}|\sigma^Tx|^2\bigg
\},\label{eq:exAV}\end{equation}
and
\[\mbox{\rm tr}(\sigma\sigma^T)=4\bigg(2+\big(1+\frac 3{\mu^2}\big
)^2\bigg),\quad |\sigma^Tx|^2\leq 4\big(1+\frac 3{\mu^2}\big)^2|x
|^2,\]
 so that, setting 
\[c_V:=4+2(\beta -1)\big(1+\frac 3{\mu^2}\big)^2,\]
\[\mathbb{A}V(x)\leq\beta |x|^{\beta -2}\bigg\{b\cdot x+c_V\bigg\}
\leq 0,\quad\forall x\in\overline {{\cal W}_0}-\{0\},\;|x|\leq\frac {
c_V}{|b|+1}.\]

Let $\delta$ be any positive number such that 
$0<\delta\leq\frac {c_V}{|b|+1}$, $0<\tilde{\delta }<\delta$ and $
0<\epsilon <\tilde{\delta}$. Also let $\tau^{X,\delta}$ and 
$\tau^{X,\epsilon}$ be defined as in $(\ref{eq:exdeltahit})$. In the sequel the 
superscript $X$ is omitted. By applying Ito's formula to the 
function $V$ (actually to a function in ${\cal C}^2_b(\R^3)$ that agrees 
with $V$ for $|x|\geq\epsilon$) we see that, for any solution $X$ of 
$(\ref{eq:bandSDER})$ in ${\cal W}_0$, with $\P\big(|X(0)|=\tilde{
\delta}\big)=1$, 
\begin{eqnarray*}
&&\E[V(X(\tau^{\epsilon}\wedge\tau^{\delta}\wedge t))]\\
&&=\big(\tilde{\delta}\big)^{\beta}+\E\big[\int_0^{\tau^{\epsilon}
\wedge\tau^{\delta}\wedge t}\mathbb{A}V(X(s))\,ds\big]+\E\big[\int_
0^{\tau^{\epsilon}\wedge\tau^{\delta}\wedge t}\nabla V(X(s))\cdot
\gamma (s)\,d\lambda (s)\big]\\
&&\leq\big(\tilde{\delta}\big)^{\beta}.\end{eqnarray*}
By sending $t\rightarrow\infty$, we get 
\[\E[V(X(\tau^{\epsilon}\wedge\tau^{\delta}))]\leq\big(\tilde{\delta}\big
)^{\beta},\]
which yields 
\[\P (\tau^{\epsilon}<\tau^{\delta})=\epsilon^{-\beta}\big(\E[V(X
(\tau^{\epsilon}\wedge\tau^{\delta}))]-\E\big[V(X(\tau^{\epsilon}
\wedge\tau^{\delta}))\mathbf{1}_{\tau^{\delta}<\tau^{\epsilon}}\big
]\big)\leq\big(\tilde{\delta}\big)^{\beta}\epsilon^{-\beta},\]
and hence, sending $\epsilon\rightarrow 0$, 
\begin{equation}\P (\vartheta^X<\tau^{X,\delta})=0.\label{eq:exmass1}\end{equation}
Now let $\tilde{\delta }=\delta /2$ and 
\[\tau^{\delta /2}_0:=0,\quad\vartheta_0:=\vartheta ,\quad\tau^{\delta}_
1:=\tau^{\delta},\quad\tau^{\delta /2}_1:=\inf\{t\geq\tau^{\delta}_
1:\,|X(t)|=\frac {\delta}2\},\]
\[\vartheta_1:=\inf\{t\geq\tau^{\delta /2}_1:\,X(t)=0\},\]
\[\tau^{
\delta}_2:=\inf\{t\geq\tau^{\delta /2}_1:\,|X(t)|=\delta \},\quad\mbox{\rm $
\tau^{\delta /2}_2:=\inf\{t\geq\tau^{\delta}_2:\,|X(t)|=\frac {\delta}
2\},$\quad etc.}\]
Then 
\[\big\{\vartheta <\infty\big\}=\bigcup_{l=0}^{\infty}\big\{\tau^{\delta /2}_l<\vartheta_l<\tau^{\delta}_{l+1}\big\}.\]
However, for $l\geq 0$, $X^l:=X(\tau^{\delta /2}_l+\cdot )$ is a solution of  
$(\ref{eq:bandSDER})$ in ${\cal W}_0$, with $\P\big(|X^l(0)|=\frac {
\delta}2\big)=1$, and 
$\vartheta_l=\vartheta^{X^l}$,  $\tau^{\delta}_{l+1}=\tau^{X^l,\delta}$, therefore, by $
(\ref{eq:exmass1})$, 
\[\P\big(\tau^{\delta /2}_l<\vartheta_l<\tau^{\delta}_{l+1}\big)=
\P\big(\vartheta^{X^l}<\tau^{X^l,\delta}\big)=0,\qquad l\geq 0,\]
and we can conclude that 
\[\P\big(\vartheta <\infty\big)=0.\]
Now let $X$ be any solution of $(\ref{eq:bandSDER})$ in 
${\cal W}_0$ with initial distribution ${\cal L}(X(0))$ supported in $\overline {
{\cal W}_0}-\{0\}$, and let 
$\delta$ be small enough that $B_{\delta}$ has empty intersection with 
the support of ${\cal L}(X(0))$. 
$X^{\delta /2}:=X(\tau^{\delta /2}+\cdot )$ is a solution of 
$(\ref{eq:bandSDER})$ in ${\cal W}_0$, with $\P\big(|X^{\delta /2}
(0)|=\frac {\delta}2\big)=1$, and 
$\vartheta^X=\vartheta^{X^{\delta /2}}$, therefore 
\[\P\big(\vartheta^X<\infty\big)=\P\big(\vartheta^{X^{\delta /2}}
<\infty\big)=0.\]

Exactly the same argument shows $(\ref{eq:ex0nothit})$.
\vskip.3in

\subsection{Starting at the origin, the origin is 
immediately left}\label{subsecexleave0}
Let $X$ be a solution of (\ref{eq:bandSDER}) in ${\cal W}_0$, that is 
$X$ satisfies $(\ref{eq:bandSDER})$ with ${\cal W}$ replaced by $
{\cal W}_0$, 
\begin{equation}\partial_4{\cal W}_0:=\partial B_4\cap {\cal W},\quad 
g^4(x):=n^4(x)=-\frac x{|x|},\label{eq:exW0}\end{equation}
and $G(x)$ is given by $(\ref{eq:bandrefcone})$, where now $h$ ranges in 
$ $$\{1,2,3,4\}$, and $(\ref{eq:bandrefcone0})$. 

Suppose $X(0)=0$ almost surely. 
Here $X$ is fixed, so we will omit the superscript $X$ in 
$\tau^{X,\delta}$. 
Let 
\[e:=\frac 1{\sqrt 3}\left[\begin{array}{c}
1\\
1\\
1\end{array}
\right],\quad f(x)=\frac 12(e\cdot x)^2.\]
By Ito's formula,  
\[\E[f(X(t\wedge\tau^{\delta}))]=\E\bigg[\int_0^{t\wedge\tau^{\delta}}\mathbb{
A}f(X(s))ds+\int_0^{t\wedge\tau^{\delta}}(e\cdot X(s))(e\cdot\gamma 
(s))\,d\lambda (s)\bigg].\]
Now, for every $x\in\partial {\cal W}_0,$ $|x|\leq\delta <4,$ 
\[e\cdot x\geq 0\mbox{\rm \ and }e\cdot g\geq 0,\quad\forall g\in 
G(x).\]
In addition 
\[\lim_{x\in\overline {{\cal W}_0}-\{0\},\,x\rightarrow 0}\mathbb{
A}f(x)=\lim_{x\in\overline {{\cal W}_0}-\{0\},\,x\rightarrow 0}\big
((e\cdot x)(e\cdot b)+\frac 12e^T(\sigma\sigma^T)e\big)=\frac 12|
\sigma^Te|^2.\]
Therefore, for $\delta$ sufficiently small, 
\[\frac 12\delta^2\geq\E[f(X(t\wedge\tau_{\delta}))]\geq\frac 14|
\sigma^Te|^2\E[t\wedge\tau^{\delta}],\]
and, by taking the limit as $t\rightarrow\infty$, 
\[\E[\tau^{\delta}]\leq\frac {2\delta^2}{|\sigma^Te|^2},\]
which of course implies also $(\ref{eq:exdeltahit})$. Exactly the same argument applies to 
an arbitrary solution of $(\ref{eq:bandSDER})$ in ${\cal W}$ 
starting at the origin. 
\vskip.3in

\subsection{If any two solutions starting at the origin have the same hitting 
distributions on $\partial B_{\delta}$, uniqueness holds for solutions 
starting at the origin}\label{subsecexhitdist}
Let $X$ be a strong Markov solution of 
$(\ref{eq:bandSDER})$ starting at the origin. Recall that 
\[\tau^{X,\delta}:=\inf\{t\geq 0:\,|X(t)|=\delta \}.\]
Since there is no risk of confusion, the superscript $X$ in 
$\tau^{X,\delta}$ will be dropped. 

Let $f$ be a continuous function on $\overline {{\cal W}_0}$ vanishing in a neighborhood 
of the origin and fix $\delta >0$ sufficiently small that $f$ 
vanishes on $B_{\delta}$. For any $\eta >0$, any solution $X$ of 
$(\ref{eq:bandSDER})$ and each $n\in\N$, we have 
\begin{eqnarray*}
&&\mathbb E\left[\int_0^{\infty}e^{-\eta t}f(X(t))dt\right]\\
&&=\mathbb 
E\left[\int_{\tau^{\delta /n}}^{\infty}e^{-\eta t}f(X(t))dt\right
]\mathbb E\left[e^{-\eta\tau^{\delta /n}}\int_{\tau^{\delta /n}}^{
\infty}e^{-\eta (t-\tau^{\delta /n})}f(X(t))\,dt\right].
\end{eqnarray*}
$\lim_{n\rightarrow\infty}\,e^{-\eta\tau^{\delta /n}}=1$ almost surely by Step 3, therefore 
\begin{eqnarray}
\label{eq:exappres}&&\mathbb E\left[\int_0^{\infty}e^{-\eta t}f(X
(t))dt\right]\\
&&=\lim_{n\rightarrow\infty}\mathbb E\left[\int_{\tau^{\delta /n}}^{
\infty}e^{-\eta (t-\tau^{\delta /n})}f(X(t))\,dt\right]=\lim_{n\rightarrow
\infty}\mathbb E\left[\int_0^{\infty}e^{-\eta t}f(X(\tau^{\delta 
/n}+t))\,dt\right]\non\end{eqnarray}
On the other hand the time-shifted process $X(\tau^{\delta /n}+\cdot 
)$ is a solution of 
$(\ref{eq:bandSDER})$ with initial distribution ${\cal L}(X(\tau^{
\delta /n}))$, 
hence, by Step 2, its distribution is uniquely 
determined by ${\cal L}\big(X(\tau^{\delta /n})\big)$. But, by assumption, $
{\cal L}\big(X(\tau^{\delta /n})\big)$ 
is the same for all solutions $X$ starting at $0$. 
Therefore, for all $n$, $\mathbb E\left[\int_0^{\infty}e^{-\eta t}
f(X(\tau^{\delta /n}+t))\,dt\right]$ is 
independent of the specific solution $X$, and hence so is 
the first line of $(\ref{eq:exappres})$.  

Then, by the uniqueness of the Laplace transform and the 
continuity of $X$, it follows 
that, for all $t\geq 0$, $\mathbb E[f(X(t))]$ is  the same for all 
solutions $X$ starting at the origin, for each function $f$ vanishing in a 
neighborhood of the origin.  
Since ${\mathbf {1}}_{\overline {{\cal W}_0}-\{0\}}$ can be approximated pointwise by a 
uniformly bounded sequence of continuous functions 
vanishing in a neighborhood of the origin, we see that, 
for each continuous $f$,  
\[\E[f(X(t))]=\E[f(X(t)){\mathbf {1}}_{\overline {{\cal W}_0}-\{0\}}(X(t))]+
f(0)\big(1-\E[{\mathbf {1}}_{\overline {{\cal W}_0}-\{0\}}(X(t))]\big)\]
is the same for all solutions $X$ starting at the origin, 
that is all solutions starting at the origin have the same 
one-dimensional distributions. By Step 2 the same 
obviously holds starting at every $x\in\overline {{\cal W}_0}-\{0
\}$. By a 
standard conditioning argument, this implies that any 
two solutions with the same initial distribution have the 
same one-dimensional distributions, and in turn this 
yields uniqueness in distribution for $(\ref{eq:bandSDER})$ 
(see Theorem \ref{th:equiv} and Lemma 3.24 of \cite{CK24aap}, 
recalling that we are considering only strong Markov 
solutions). 
\vskip.3in

\subsection{The ergodic argument and conclusion of the 
proof of Theorem \ref{th:exuniq}}\label{subsecexergodic}

For each $n\in\N$, the inhomogeneous Markov chain 
$\big\{X_k\big\}_{0\leq k\leq n}:=\big\{X(\tau^{n-k})\big\}_{0\leq 
k\leq n}$, where $\tau^{n-k}:=\tau^{X,2^{-2(n-k)}\delta}$ 
is the first time $X$ hits $\partial B_{2^{-2(n-k)}\delta}$, has transition kernels 
\begin{equation}\begin{array}{c}
Q_{n-k}(x,B):=\P\big(X(\tau^{n-k-1})\in B\big|X(\tau^{n-k})=x\big
),\\
\hfill\\
B\in {\cal B}\big(\partial B_{2^{-2(n-k-1)}\delta}\big),\;x\in\partial 
B_{2^{-2(n-k)}\delta},\;0\leq k\leq n-1.
\end{array}\label{eq:extransker}\end{equation}
Then, denoting by $Q_{n-k}$ the integral operators 
corresponding to the kernels $(\ref{eq:extransker})$, for every $
f\in {\cal C}(\partial B_{\delta})$, 
\begin{equation}\E\big[f(X(\tau^{\delta}))\big]=\E\big[f(X_n)\big
]=\int_{\partial B_{2^{-2n}\delta}}Q_n\cdots Q_1f(x)\,d\nu_n(x),\label{eq:exbasid}\end{equation}
where $\nu_n;={\cal L}(X_0)={\cal L}(X(\tau^n)$). As discussed in Step 5) of 
Section \ref{subsecexoutline}, for another solution 
$\tilde {X}$ of $(\ref{eq:bandSDER})$ starting at the origin, the 
transition kernels of the corresponding Markov chain $\big\{\tilde {
X}_k\big\}_{0\leq k\leq n}$ are 
the same as for $X$, so that 
\[\E\big[f(\tilde {X}(\tilde{\tau}^{\delta}))\big]=\E\big[f(\tilde {
X}_n)\big]=\int_{\partial B_{2^{-2n}\delta}}Q_n\cdots Q_1f(x)\,d\tilde{
\nu}_n(x),\]
where $\tilde{\nu}_n;={\cal L}(\tilde {X}_0)={\cal L}(\tilde {X}(
\tilde{\tau}^n)$). 

By the reverse ergodic theorem of 
\cite{CK24aap} (recalled below as Theorem \ref{th:ergodic}), the limit 
\[\lim_{n\rightarrow\infty}\int_{\partial B_{2^{-2n}\delta}}Q_n\cdots 
Q_1f(x)\,d\nu_n(x)\]
exists and is independent of the sequence of probability 
measures $\big\{\nu_n\big\}_{n\geq 1}$, therefore 
\[\E\big[f(X(\tau^{\delta})\big]=\E\big[f(\tilde {X}(\tilde{\tau}^{
\delta})\big],\]
for every $f\in {\cal C}(\partial B_{\delta})$, that is the hitting distributions of 
$X$ and $\tilde {X}$ coincide. 
\vskip.3in

\section{Uniqueness for semimartingale ORBM in a 
piecewise ${\cal C}^2$ cone} 
\label{sectionpiecesmcone}

\setcounter{equation}{0}

The study of semimartingale obliquely reflecting Brownian 
motion (ORBM) in a piecewise smooth cone, i.e. the intersection 
of a finite number of smooth cones,  
cannot be reduced to the study of semimartingale ORBM in a piecewise 
smooth domain.  
Therefore the results of \cite{DI93} or 
\cite{DW95} do not apply. Moreover, the arguments used 
in these papers do not carry over to piecewise smooth 
cones. 

\cite{KW91} has studied semimartingale ORBM in a smooth cone ${\cal W}$, with 
radially constant, smooth direction of reflection $g$. The 
results of this paper use the smoothness of the cone in two ways: 

\noindent- to find a function $\Phi$ that solves 
\begin{equation}\left\{\begin{array}{ll}
\Delta\Phi =0,&\mbox{\rm in }{\cal W},\\
\nabla\Phi\cdot g=0,&\mbox{\rm on }\partial {\cal W}-\{0\},\end{array}
\right.\label{eq:KW}\end{equation}
where $0$ is the vertex of the cone; the function $\Phi$ 
determines whether the vertex is reached;

\noindent- to prove the support theorem (Theorem 3.1 of 
\cite{KW91}), which allows to use 
the Krein-Rutman theorem for strongly positive operators 
to uniquely characterize the ORBM. 

Here we seek to uniquely characterize the 
semimartingale ORBM in a 
piecewise smooth cone ${\cal W}$. 
The main difference between this work and 
\cite{KW91} is that the present work relies not on the 
Krein-Rutman theorem, but on a new reverse, ergodic 
theorem for killed Markov chains proved in \cite{CK24aap} 
and recalled below as Theorem \ref{th:ergodic}. The 
assumptions of Theorem \ref{th:ergodic} have a clear 
probabilistic meaning and may be verified in a wider 
range of situations: piecewise smooth cones - the case 
of this paper; directions of reflection that are not 
radially constant; variable coefficients of the reflecting 
process; domains that are not cones but can be locally 
approximated by cones (see \cite{CK24aap} and 
\cite{CK24spa}), etc.. Moreover verifying the assumptions of Theorem 
\ref{th:ergodic} does not 
require the delicate oscillation estimates needed to apply 
the Krein-Rutman theorem (see Theorem 3.2 of 
\cite{KW91}).

As in \cite{KW91}, here too we employ a function (or two 
functions) to analyze the time till the vertex is 
reached, but we use the observation that, for our 
purposes, the functions do not need to satisfy the 
equalities in $(\ref{eq:KW})$, but only corresponding 
inequalities (Condition \ref{auxfunc} below). 
Thus one does not need to solve the boundary 
value problem $(\ref{eq:KW})$, but only to find functions 
that satisfy Condition \ref{auxfunc},  
and one can avoid requiring 
smoothness of the cone (as an example, see the function 
in Section \ref{subsecex0nothit}). 
On the other hand, while 
\cite{KW91} gives a precise classification of the possible 
behaviours of the time till the vertex is reached, here 
we only obtain upper and lower bounds on 
$\P (\tau^{X,\delta}<\vartheta^X)$ (where $X$ is a semimartingale ORBM starting at 
$x$, $0<|x|<\delta$, $\tau^{X,\delta}:=\inf\{t\geq 0:\;|X(t)|=\delta 
\}$, 
$\vartheta^X:=\inf\{t\geq 0:\;X(t)=0\}$). However these bounds 
are enough to satisfy the first assumption of Theorem 
\ref{th:ergodic}. In particular one can give a 
sufficient condition for $\P (\tau^{X,\delta}<\vartheta^X)=1$, which implies 
that, starting at $x\neq 0$, the origin is never 
reached (this is the case of the example of Section 
\ref{sectionex}: see Section \ref{subsecex0nothit}). 

The second assumption of Theorem \ref{th:ergodic} 
replaces the assumption of the Krein-Rutman theorem 
that the operator (in applications to ORBM, the transition 
operator of a killed Markov chain) 
be a compact, strongly positive 
operator. This second assumption is 
verified by a coupling lemma proved in 
\cite{CK18} (Lemma 5.3 of \cite{CK18}) and by 
a scaling limit and a nontrivial extension of the support theorem of \cite{KW91} to 
piecewise ${\cal C}^2$ cones (Lemma \ref{th:scaling} below). 
A careful inspection shows that the 
proof of the support theorem of \cite{KW91} carries 
over if one can prove that the semimartingale ORBM obtained 
from the scaling limit spends zero time on $\partial {\cal W}-\{0
\}$ 
before hitting the origin. 
This follows from a result 
of \cite{KR17} (Proposition 2.12 of \cite{KR17}) if ${\cal W}$ is piecewise $
{\cal C}^2$, and this is why ${\cal C}^1$ 
smoothness of the faces of ${\cal W}$ is not sufficient. 

The problem, the assumptions and the main result of the 
paper (Theorem \ref{th:uniq}) are formulated precisely in Section 
\ref{subsecassum} below. At the end of 
Section \ref{subsecassum} the outline of the 
proof of Theorem \ref{th:uniq} is discussed and compared 
with that of Theorem \ref{th:exuniq}. 
The proof is then carried out in Sections 
\ref{subseclocaliz} to \ref{subsecscale}. Section \ref{subsecCMP} 
contains a preliminary result and the Appendix 
contains the proof of a technical result.

\subsection{Formulation of the problem and main 
result}\label{subsecassum}
Let ${\cal W}\subseteq\R^d$ be a piecewise ${\cal C}^2$ cone with 
vertex at the origin i.e. 
\begin{equation}{\cal W}:=\bigcap_{j=1}^m{\cal W}_j,\quad {\cal W}_
j:=\{x=rz,\,z\in {\cal S}_j,\,r>0\},\label{eq:cone}\end{equation}
where ${\cal S}_j$ is a nonempty domain in the unit sphere $S^{d-
1}$ 
with ${\cal C}^2$ boundary. More precisely 
\begin{equation}\begin{array}{c}
{\cal S}_j=\{z\in S^{d-1}:\,\varphi_j(z)>0\},\quad\partial {\cal S}_
j=\{z\in S^{d-1}:\,\varphi_j(z)=0\},\\
\varphi_j\in {\cal C}^2(S^{d-1}),\qquad\inf_{z\in\partial {\cal S}_
j}|\nabla_{S^{d-1}}\varphi_j(z)|>0.\end{array}
\label{eq:domfunc}\end{equation}
Clearly 
\begin{equation}{\cal W}=\{x=rz,\,z\in {\cal S},\,r>0\},\;\;\;\mbox{\rm where }
{\cal S}:=\bigcap_{j=1}^m{\cal S}_j.\label{eq:cone2}\end{equation}
Let $g^j$ be a radially constant direction of 
reflection on the face 
\[\partial_j{\cal W}:=\partial {\cal W}_j\cap\bigcap_{i\neq j}{\cal W}_
i,\]
$ $that is $g^j$ is a unit vector field defined on $\partial {\cal W}_
j-\{0\}$ and, for 
$x=rz$, $z\in\partial {\cal S}_j$, $r>0$,
\begin{equation}g^j(x)=g^j(z).\label{eq:radhg}\end{equation}
A semimartingale ORBM can be defined as a solution of the 
following stochastic differential equation with reflection: 
\begin{eqnarray}
&X(t)=X(0)+b\,t+\sigma\,W(t)+\int_0^t\gamma (s)\,d\lambda (s),\quad 
t\geq 0,&\non\\
&X(t)\in\overline {{\cal W}},\quad\lambda (t)=\int_0^t{\bf 1}_{\partial 
{\cal W}}(X(s))d\lambda (s),\quad t\geq 0,&\label{eq:SDER}\\
&\quad\gamma (t)\in G(X(t)),\quad |\gamma (t)|=1,\quad d\lambda -
a.e.,\quad t\geq 0,&\non\end{eqnarray}
where, for $x\in\partial {\cal W}-0$, $G(x)$ is the cone of possible 
directions of reflection at $x$:
\begin{equation}G(x):=\{g:\,g=\sum_{j\in {\cal J}(x)}u_jg^j(x),\;
u_j\geq 0\},\quad {\cal J}(x):=\{j:\,x\in\overline {\partial_j{\cal W}}
\},\label{eq:refcone}\end{equation}
and $G(0)$ is the cone of possible directions of reflection 
at $0$, that is the closed convex cone generated by 
\begin{equation}\big\{g^j(x),\;x\in\ \overline {\partial_j{\cal W}}
-\{0\},\;j=1,\cdots ,m\big\}.\label{eq:refcone0}\end{equation}

\begin{definition}\label{def:SDER}

A continuous process $X$ is a solution of $(\ref{eq:SDER})$ 
if there exist a standard Brownian motion $W$, a continuous, 
nondecreasing process $\lambda$ and a process 
$\gamma$ with measurable paths, all defined on the same probability space as $
X$, 
such that $W(t+\cdot )-W(t)$ is independent of 
${\cal F}_t^{X,W,\lambda ,\gamma}$ for all $t\geq 0$ and $(\ref{eq:SDER})$ 
is satisfied a.s.. 

{\em Weak uniqueness\/} or {\em uniqueness in distribution\/} holds if any 
two solutions of $(\ref{eq:SDER})$ with
the same initial distribution have the same distribution on 
$C_{\overline {{\cal W}}}[0,\infty )$. 

A stochastic process $\tilde {X}$ (for example a solution of an 
appropriate martingale problem or submartingale problem) 
is a weak solution of $(\ref{eq:SDER})$ if there is a 
solution $X$ of (\ref{eq:SDER}) such that $\tilde {X}$ and $X$ have the 
same distribution.
\end{definition}

\begin{remark}\label{re:representation}
There is not a unique representation 
$(\ref{eq:cone})$-$(\ref{eq:domfunc})$ of ${\cal W}$, 
as it would be for a polyhedral cone, because each 
function $\varphi_j$ can be extended from $\overline {{\cal S}}$ do $
S^{d-1}$ in infinitely 
many ways. However, we only need that there exists 
at least one representation and a set of directions of 
reflection $\{g^j\}_{j=1,\cdots ,m}$ that satisfy Conditions \ref{G} 
and \ref{auxfunc} below. 
\end{remark}
\vskip.2in

The goal of this section is to provide conditions on ${\cal W}$, $
G$ 
and $\sigma$ under which uniqueness in distribution holds for $(\ref{eq:SDER})$. 
Since the focus is on uniqueness, existence of solutions to 
$(\ref{eq:SDER})$ for every initial distribution 
$\nu\in {\cal P}(\overline {{\cal W}})$ will be assumed throughout. 

For $x\in\partial {\cal W}_j-\{0\}$, $n^j(x)$ will denote the unit inward 
normal to $\overline {{\cal W}_j}$ and, for each $x\in\partial {\cal W}
-\{0\}$, 
$N(x)$ will denote the closed, convex cone generated by $\big\{n^
j(x)\big\}_{j\in {\cal J}(x)}$. 
Note that, for $x\in\partial {\cal W}_j-\{0\}$, $x=rz$, $z\in\partial 
{\cal S}_j$, $r>0$, 
\begin{equation}n^j(x)=n^j(z),\label{eq:radhn}\end{equation}
and $n^j(z)$ concides with the inward unit normal to $\overline {
{\cal S}_j}$ in 
$S^{d-1}$, in particular $n^j(z)$ lies in the tangent space to 
$S^{d-1}$ at $z.$ Also note that, for $x\in\partial {\cal W}-\{0\}$, $
x=rz$, $z\in\partial {\cal S}$, $r>0$, 
\begin{equation}{\cal J}(x)={\cal J}(z):=\{j:\,z\in\partial {\cal S}_
j\}.\label{eq:J}\end{equation}
For each $\zeta\in\partial {\cal S}$, if $j\notin {\cal J}(\zeta 
)$, $\zeta$ has positive distance from 
$\partial {\cal S}_j$, hence there is $\delta_j(\zeta )>0$ such that, for 
$z\in\partial {\cal S}\cap B_{\delta_j(\zeta )}(\zeta )$, $z$ has positive distance from $
\partial {\cal S}_j$ as well, 
so that $j\notin {\cal J}(z)$, therefore, setting $\delta_{{\cal J}}
(\zeta ):=\min_{j\notin {\cal J}(\zeta )}\delta_j(\zeta )$, 
\begin{equation}{\cal J}(z)\subseteq {\cal J}(\zeta ),\quad\forall 
z\in\partial {\cal S}\cap B_{\delta_{{\cal J}}(\zeta )}(\zeta ).\label{eq:deltaJ}\end{equation}

The main assumptions are formulated in the following 
two conditions. 

\begin{condition}\label{G}\hfill 
\begin{itemize}

\item[(i)]For $j=1,\ldots ,m$, $g^j:S^{d-1}\rightarrow\R^d$ is a Lipschitz continuous  
vector field of unit length on $\partial {\cal S}_j$ such that 
\[\inf_{z\in\partial {\cal S}_j}g^j(z)\cdot n^j(z)>0.\]
\item[(ii)]For every $z\in\partial {\cal S}$, the vectors $\{n^j(
z)\}_{j\in {\cal J}(z)}$ are 
linearly independent. In particular $|{\cal J}(z)|\leq d-1$. 

\item[(iii)]For $z\in\partial {\cal S}$, ${\cal J}(z)=\{h_1,\ldots 
,h_k\}$, the matrix of elements 
\[\frac {|g^{h_i}(z)\cdot n^{h_j}(z)|}{g^{h_i}(z)\cdot n^{h_i}(z)}
-\delta_{ij},\quad i,j=1,...,k,\]
has spectral radius strictly less than $1$.
\item[(iv)]Let 
\begin{equation}N(0):=\{n\in\R^d:\,n\cdot x\geq 0,\,\forall x\in\overline {
{\cal W}}\},\label{eq:N(0)}\end{equation}
and suppose that $\overset{\circ}N(0)\neq\emptyset$. 

There exists a unit vector $e\in N(0)$ such that 
\[e\cdot g>0,\quad\forall g\in G(0),\,|g|=1.\]
Without loss of generality, one can suppose 
$e\in\overset{\circ}N(0)$. 
\end{itemize}
\end{condition}
\vskip.2in

\begin{remark}\label{re:e(z)}
For each $z\in\partial {\cal S}$, let 
\[N(z):=\{n:\,n=\sum_{j\in {\cal J}(z)}u_jn^j(z),\;u_j\geq 0\}.\]
As in the Remark at page 160 of \cite{DI91}, perturbing 
the matrix  of Condition \ref{G} (iii) by a matrix with 
small positive entries and using the Perron-Frobenius 
theorem, we see that for each $z\in\partial {\cal S}$ there is a unit 
vector $e(z)\in N(z)$ and $\delta (z)>0$, such that 
\begin{equation}e(z)\cdot g\geq c_e>0,\qquad\forall g\in G(\zeta 
),\,|g|=1,\;\forall\zeta\in\partial {\cal S},\,|\zeta -z|\leq\delta 
(z),\label{eq:e(z)}\end{equation}
where $c_e$ is a positive constant independent of $z$.
\end{remark}
\vskip.2in

In Condition \ref{G}, assumptions (i) to (iii) guarantee that, 
at every point of the boundary other than the origin the 
assumptions of \cite{DI93} are satisfied. In particular 
$(\ref{eq:e(z)})$ ensures that, if a solution of 
$(\ref{eq:SDER})$ reaches $\partial {\cal W}-\{0\}$, it leaves the boundary 
immediately.  Assumption (iv) 
ensures that the same holds at the origin. 
\vskip.2in

Denote by $\mathbb{A}$ the operator of the form 
\begin{equation}\mathbb{A}f(x):=b\cdot\nabla f(x)+\frac 12\mbox{\rm tr}
((\sigma\sigma^T)D^2f(x)).\label{eq:A}\end{equation}
$\sigma$ is assumed to be nonsingular.
\vskip.2in

\begin{condition}\label{auxfunc}\hfill 

For some $\delta_{{\cal W}}>0$, either of the following conditions is satisfied
\begin{itemize}
\item[(i)]There exists a function $V\in {\cal C}^2(\overline {{\cal W}}
-\{0\})$ such 
that 
\begin{equation}\lim_{x\in\overline {{\cal W}},\,x\rightarrow 0}V
(x)=\infty ,\label{eq:V-0}\end{equation}
\begin{equation}\begin{array}{c}
\nabla V(x)\cdot g\leq 0,\quad\forall\,g\in G(x),\,
x\in\big(\partial {\cal W}-\{0\}\big)\cap\overline {B_{\delta_{{\cal W}}}
(0)},\\
\mathbb{A}V(x)\leq 0,\quad\forall\,x\in\big(\overline {
{\cal W}}-\{0\}\big)\cap\overline {B_{\delta_{{\cal W}}}(0)}.
\end{array}\label{eq:nablaAV-}\end{equation}
\item[(ii)]There exist two functions 
$V_{+},V_{-}\in {\cal C}^2(\overline {{\cal W}}-\{0\})$ such that 
\begin{equation}\begin{array}{c}
V_{+}(x)>0,\quad V_{-}(x)>0,\quad\mbox{\rm for }x\in\big(\overline {
{\cal W}}-\{0\}\big)\cap\overline {B_{\delta_{{\cal W}}}(0)},\\
\lim_{x\in\overline {{\cal W}},\,x\rightarrow 0}V_{+}(x)=\lim_{x\in\overline {
{\cal W}},\,x\rightarrow 0}V_{-}(x)=0,\end{array}
\label{eq:V+0}\end{equation}
\begin{equation}\inf_{0<\delta\leq\delta_{{\cal W}}}\frac {\inf_{
|x|=\delta}V_{+}(x)}{\sup_{|x|=\delta}V_{-}(x)}>0,\quad\inf_{0<\delta
\leq\delta_{{\cal W}}}\frac {\inf_{|x|=\delta}V_{-}(x)}{\sup_{|x|
=\delta}V_{+}(x)},>0\label{eq:V+unif}\end{equation}
\begin{equation}\begin{array}{c}
\nabla V_{+}(x)\cdot g\geq 0,\quad\nabla V_{-}(x)\cdot g\leq 0,\quad
\forall\,g\in G(x),\,x\in\big(\partial {\cal W}-\{0\}\big)\cap\overline {
B_{\delta_{{\cal W}}}(0)}\\
\hfill \\
\mathbb{A}V_{+}(x)\geq 0,\quad\mathbb{A}V_{-}(x)\leq 0,\quad\forall\,
x\in\big(\overline {{\cal W}}-\{0\}\big)\cap\overline {B_{\delta_{
{\cal W}}}(0)}.\end{array}
\label{eq:nablaAV+}\end{equation}
\end{itemize}
\end{condition}
\vskip.2in

For a solution $X$ of $(\ref{eq:SDER})$, let 
\begin{equation}\vartheta^X:=\inf\{t\geq 0:\,X(t)=0\},\label{eq:0hittime}\end{equation}
\begin{equation}\tau^{X,\delta}:=\inf\{t\geq 0:\,|X(t)|=\delta \}
,\quad\delta >0.\label{eq:deltahittime}\end{equation}
(When there is no risk of confusion, the superscript $X$ 
will be omitted.). In addition denote by $X^x$ a solution of 
$(\ref{eq:SDER})$ starting at $x$. 
Condition \ref{auxfunc} (i) is a sufficient condition 
for the following: for $\delta\leq\delta_{{\cal W}}$ and $0<|x|<\delta$, 
\begin{equation}\P (\tau^{X^x,\delta}<\vartheta^{X^x})=1.\label{eq:V-}\end{equation}
In particular this yields 
\[\P (\vartheta^{X^x}<\infty )=0,\qquad\forall\,x\neq 0\]
(see Section \ref{subsecex0nothit}). 

The meaning of Condition \ref{auxfunc} (ii) instead is the 
following. If Condition \ref{auxfunc} (ii) is satisfied, for 
$\delta\leq\delta_{{\cal W}}$, 
\begin{equation}\lim_{x\rightarrow 0}\P (\tau^{X^x,\delta}<\vartheta^{X^x})=0,\label{eq:V+}\end{equation}
but the above probability vanishes at the same rate for all 
points that have the same distance from the origin, that 
is the following uniform lower bound holds:
\[\inf_{x,\tilde {x}\in\overline {{\cal W}}:\,0<|x|=|\tilde {x}|<
\delta}\frac {\P (\tau^{X^x,\delta}<\vartheta^{X^x})}{\P (\tau^{X^{
\tilde {x}},\delta}<\vartheta^{X^{\tilde {x}}})}\geq c_0>0.\]
Of course this trivially holds also if Condition \ref{auxfunc} (i) 
is satisfied. In the context of this work, the above uniform lower bound is 
the first assumption of Theorem \ref{th:ergodic}. 

Note that $(\ref{eq:V-})$ and $(\ref{eq:V+})$ imply that 
only one between Conditions \ref{auxfunc} (i) and (ii) can 
be verified.  
\vskip.2in

The main result of this work is the following.  

\begin{theorem}\label{th:uniq}
Under Conditions \ref{G} and \ref{auxfunc}, for every 
$\nu\in {\cal P}(\overline {{\cal W}})$ uniqueness in 
distribution holds for $(\ref{eq:SDER})$. 

The solution, $X$, spends zero time at the origin, that is 
\begin{equation}\int_0^{\infty}\mathbf{1}_{\{0\}}(X(t))\,dt=0,\qquad 
a.s..\label{eq:0zerotime}\end{equation}

$(\ref{eq:0zerotime})$ implies in particular that, if $X$ starts at $
0$, 
\begin{equation}\lim_{\delta\rightarrow 0}\tau^{X,\delta}=0,\qquad 
a.s..\label{eq:deltahit}\end{equation}
\end{theorem}

\begin{remark}\label{re:nonradconst}
Assuming that the directions of reflection are radially 
constant (that is $(\ref{eq:radhg})$ holds) simplifies proofs. 
However, it can be shown, by the same arguments as in 
\cite{CK24aap}, that Theorem \ref{th:uniq} carries over to radially 
variable directions of reflection, as long as 
\[|g^j(rz)-\bar {g}^j(z)|\leq c_{{\cal W}}r,\qquad\forall z\in\partial 
{\cal S}_j,\,0<r\leq r_{{\cal W}},\]
for some unit vector field $\bar {g}^j\in {\cal C}^2(S^{d-1},\R^d
)$ of unit 
length on $\partial {\cal S}_j$ and some $c_{{\cal W}},r_{{\cal W}}
>0$.
\end{remark}
\vskip.2in

Theorem \ref{th:uniq} is a general result which includes 
Theorem \ref{th:exuniq} as a special case. The proof 
follows the outline given in Section 
\ref{subsecexoutline}, but Step 5 is more complex. 
In fact, in general, the solutions of 
$(\ref{eq:SDER})$ starting at the origin can hit the origin between the hitting 
time of $\partial B_{2^{-2(n-k)}\delta}$ and the hitting time of $
\partial B_{2^{-2(n-k-1)}\delta}$, so 
that, as mentioned in Remark \ref{re:exkilled}, in order 
to have transition kernels that are independent of the 
specific solution, one has to consider 
Markov chains killed when the origin is reached. 
However the reverse ergodic theorem of \cite{CK24aap} 
applies to killed Markov chains as well, under a suitable 
assumption which is ensured by Condition \ref{auxfunc} (see 
Section \ref{subseczerotime0}). Step 3 is also 
partly different, because, in order to carry out Step 4 
in the general case, one needs to prove also 
that any solution of $(\ref{eq:SDER})$
spends zero time at the origin, 

Each of the steps of the proof of 
Theorem \ref{th:uniq} is developed in one of the 
Sections \ref{subseclocaliz} to \ref{subsecscale}. 
While in Section \ref{sectionex} 
only key points were proved, but those were proved in 
detail, here proofs are complete but concise and often 
refer to previous works. Section 
\ref{subsecCMP} contains some preliminary 
results.

\subsection{Constrained martingale problems}\label{subsecCMP}

ORBMs (and more generally stochastic processes with 
boundary conditions) can be 
characterized also as {\em natural solutions of constrained martingale problems}. 
Characterizing the process 
as the natural solution of a constrained martingale problem 
allows to verify tightness in a much simpler way, 
whenever an approximation or convergence argument is needed. 
The present work builds on some previous results which are 
stated in terms of constrained martingale problems, 
therefore definitions and some of the results are recalled here. 
Definitions and results are formulated for the specific 
constrained martingale problem corresponding to 
$(\ref{eq:SDER})$. The reader is referred to 
\cite{Kur90}, \cite{Kur91}, \cite{CK15}, \cite{CK19} for 
a complete treatment. 
In these papers the state space is compact, but the 
definitions extend to a locally compact state space 
as well as most of the results. 

Let $ $$E_0$ be an open subset of $\R^d$, 
\begin{equation}{\cal D}:={\cal C}^2_b(\overline {E_0}),\label{eq:dom}\end{equation}
$ $and $\A$ be given by $(\ref{eq:A})$ with domain 
\begin{equation}{\cal D}(\A)={\cal D}.\label{eq:Adom}\end{equation}
To each $x\in\partial E_0$ associate a convex cone $G(x)$ in such a way that the set 
\begin{equation}\Xi :=\{(x,u)\in\partial E_0\times S^{d-1}:u\in G
(x)\},\label{eq:clset}\end{equation}
is closed and define 
\begin{equation}\mathbb{B}:{\cal D}\rightarrow\,{\cal C}(\Xi ),\quad\mathbb{
B}f(x,u):=\nabla f(x)\cdot u\label{eq:B}\end{equation}

Define ${\cal L}_{\Xi}$ to be the space of measures $\mu$ on 
$[0,\infty )\times\Xi$ such that $\mu ([0,t]\times\Xi )<\infty$ for all $
t>0$.  ${\cal L}_{\Xi}$  
is topologized so that $\mu_n\in {\cal L}_{\Xi}\rightarrow\mu\in 
{\cal L}_{\Xi}$ if and only if 
\[\int_{[0,\infty )\times\Xi}f(s,u)\mu_n(ds\times du)\rightarrow\int_{
[0,\infty )\times\Xi}f(s,u)\mu (ds\times du)\]
for all continuous $f$ with compact support in $[0,\infty )\times
\Xi$. It is 
possible to define a metric on ${\cal L}_{\Xi}$ that induces the above 
topology and makes ${\cal L}_{\Xi}$ into a complete, separable metric 
space. Also define ${\cal L}_{S^{d-1}}$ in the same way. 

An ${\cal L}_{S^{d-1}}$-valued random variable $\Lambda_1$ is adapted to a filtration $
\{{\cal G}_t\}$ if 
\[\Lambda_1([0,\cdot ]\times B)\mbox{\rm \ is }\{{\cal G}_t\}-\mbox{\rm adapted}
,\quad\forall B\in {\cal B}(S^{d-1}).\]
An adapted ${\cal L}_{\Xi}$-valued random variable is defined 
analogously. 

\begin{definition}\label{def:cmp}
Let $\mathbb{A}$, $\Xi$ and $\mathbb{B}$ be as in $(\ref{eq:A})$, 
$(\ref{eq:Adom})$, $(\ref{eq:clset})$ and $(\ref{eq:B})$. 
A process $X$ in $D_{\overline {E_0}}[0,\infty )$ is a solution of the 
{\em constrained martingale problem\/} for $(\mathbb{A},E_0,\mathbb{
B},\Xi )$ if 
there exists a random measure $\Lambda$ with values in ${\cal L}_{
\Xi}$ and a filtration 
$\{{\cal F}_t\}$ such that $X$ and $\Lambda$ are $\{{\cal F}_t\}$-adapted and for each 
$f\in {\cal D}$, 
\begin{equation}f(X(t))-f(X(0))-\int_0^t\mathbb{A}f(X(s))ds-\int_{
[0,t]\times\Xi}\mathbb{B}f(x,u)\Lambda (ds\times dx\times du)\label{eq:cmp}\end{equation}
is a $\{{\cal F}_t\}$-local martingale. 
\end{definition}

A {\em natural solution\/} of a constrained martingale problem 
is a solution obtained by time-changing 
a solution of the corresponding {\em controlled martingale }
{\em problem}, which is a ''slowed down'' version of the 
constrained martingale problem. 

\begin{definition}\label{def:clmp}
Let $\mathbb{A}$, $\Xi$ and $\mathbb{B}$ be as in $(\ref{eq:A})$, 
$(\ref{eq:Adom})$, $(\ref{eq:clset})$ and $(\ref{eq:B})$. 
$(Y,\lambda_0,\Lambda_1)$ is a solution of the {\em controlled martingale problem }
for $(\mathbb{A},E_0,\mathbb{B},\Xi )$, if $Y$ is a process in 
${\cal D}_{\overline {E_0}}[0,\infty )$, $\lambda_0$ is nonnegative and nondecreasing, 
$\Lambda_1$ is a random measure with values in ${\cal L}_{S^{d-1}}$ such that 
\begin{equation}\lambda_1(t):=\Lambda_1([0,t]\times S^{d-1})=\int_{
[0,t]\times S^{d-1}}{\bf 1}_{\Xi}(Y(s),u)\Lambda_1(ds\times du),\label{Lam1supp}\end{equation}
\[\lambda_0(t)+\lambda_1(t)=t,\]
and there exists a filtration $\{{\cal G}_t\}$ such that $Y$, $\lambda_
0$, and $\Lambda_1$ are 
$\{{\cal G}_t\}$-adapted and 
\begin{equation}f(Y(t))-f(Y(0))-\int_0^t\mathbb{A}f(Y(s))d\lambda_
0(s)-\int_{[0,t]\times S^{d-1}}\mathbb{B}f(Y(s),u)\Lambda_1(ds\times 
du)\label{eq:clmgp}\end{equation}
is a $\{{\cal G}_t\}$-martingale for all $f\in {\cal D}$. 
\end{definition}

For every nondecreasing path $\lambda_0\in {\cal D}_{[0,\infty )}
[0,\infty )$, define 
\begin{equation}\lambda_0^{-1}(t)=\inf\{s:\lambda_0(s)>t\},\quad 
t\geq 0.\label{eq:tinv}\end{equation}
\begin{definition}\label{def:nat}
Let $\mathbb{A}$, $\Xi$ and $\mathbb{B}$ be as in $(\ref{eq:A})$, 
$(\ref{eq:Adom})$, $(\ref{eq:clset})$ and $(\ref{eq:B})$. 
A solution, $X$, of the constrained martingale problem 
for $(\mathbb{A},E_0,\mathbb{B},\Xi )$ is called {\em natural\/} if, for some 
solution $(Y,\lambda_0,\Lambda_1)$ of the controlled martingale problem 
for $(\mathbb{A},E_0,\mathbb{B},\Xi )$ with filtration $\{{\cal G}_
t\}$, 
\begin{eqnarray}
\label{eq:tchange}&&\qquad\qquad\qquad X(t)=Y(\lambda_0^{-1}(t)),
\quad {\cal F}_t={\cal G}_{\lambda_0^{-1}(t)},\non\\
&&\Lambda ([0,t]\times C)=\int_{[0,\lambda_0^{-1}(t)]\times S^{d-
1}}{\bf 1}_C(Y(s),u)\Lambda_1(ds\times du),\quad C\in {\cal B}(\Xi 
),\end{eqnarray}
Uniqueness holds for natural solutions of 
the constrained martingale problem for 
$(\mathbb{A},E_0,\mathbb{B},\Xi )$ if any two natural 
solutions 
with the same initial distribution have the same 
distribution on $D_{\overline {E_0}}[0,\infty )$. 
\end{definition}

The two characterizations of a semimartingale ORBM as a 
solution of $(\ref{eq:SDER})$ and as a natural solution of 
the constrained martingale problem for 
$(\mathbb{A},{\cal W},\mathbb{B},\Xi )$, with the cone $G(x)$ given 
by $(\ref{eq:refcone})$ and $(\ref{eq:refcone0})$, are equivalent. 
In particular, taking into account Lemma 
\ref{th:lambda0pos} below, if $X$ is a natural solution 
of the constrained martingale problem 
for $(\mathbb{A},E_0,\mathbb{B},\Xi )$, by $(\ref{eq:tchange})$ 
and $(\ref{Lam1supp})$ the corresponding $\Lambda$ satisfies 
\begin{eqnarray*}
&&\int_{[0,t]\times\Xi}h(x,u)\,\Lambda (dr\times dx\times du)=\int_{
[0,\lambda_0^{-1}(t)]\times S^{d-1}}h(Y(r),u)\,\Lambda_1(dr\times 
du)\\
&&=\int_{[0,t]\times\Xi}h(X(r),u)\,\Lambda (dr\times dx\times du)\end{eqnarray*}
Therefore $\Lambda$ is concentrated on the set 
\[\big\{(r,x,u)\in [0,t]\times\Xi :\,x=X(r)\big\}\\
\]
\[=\big\{(r,x,u):\,0\leq r\leq t,\,x=X(r)\in\partial {\cal W},\,u
\in G(X(r)),\,|u|=1\big\},\]
which corresponds to the conditions in the second and third 
lines of $(\ref{eq:SDER})$ (except $X(t)\in\overline {{\cal W}}$, which holds by 
definition for any solution). 

Theorem \ref{th:equiv} below formulates the equivalence 
precisely. 

\begin{lemma}\label{th:lambda0pos}
For every solution $(Y,\lambda_0,\Lambda_1)$ of the controlled martingale problem 
for $(\mathbb{A},{\cal W},\mathbb{B},\Xi )$, $\lambda_0(t)>0$ for all $
t>0$, 
almost surely. Moreover $\lambda_0$ is strictly increasing almost 
surely. 
\end{lemma}

\begin{proof}
The proof that $\lambda_0(t)>0$ for all $t>0$ is analogous to that of Lemma 6.8 of 
\cite{CK19} and Lemma 3.1 of \cite{DW95} and is based on 
Remark \ref{re:e(z)} and 
Condition \ref{G} (iv). The second statement follows by Lemma 3.4 of 
\cite{CK19} (note that Lemma 3.4 of \cite{CK19} 
holds in locally compact spaces as well).  
\end{proof}

\begin{theorem}\label{th:equiv}
Let $\Xi$ be as in $(\ref{eq:clset})$ with $G(x)$, $x\in\partial 
{\cal W}$, given 
by $(\ref{eq:refcone})$ and $(\ref{eq:refcone0})$, 
$\mathbb{A}$, $\mathbb{B}$ be as in $(\ref{eq:A})$, 
$(\ref{eq:B})$, and assume Conditions \ref{G} 
and \ref{auxfunc}.\hfill\break 
Every solution of $(\ref{eq:SDER})$ is a 
natural solution of the constrained martingale problem 
for $(\mathbb{A},{\cal W},\mathbb{B},\Xi )$.\break
Every natural solution of the constrained martingale problem 
for $(\mathbb{A},{\cal W},\mathbb{B},\Xi )$ is a weak solution of 
$(\ref{eq:SDER})$. 
\end{theorem}

\begin{proof}
The proof relies on Lemma \ref{th:lambda0pos} and is the same as 
for Theorem 6.12 of \cite{CK19}. 
\end{proof}

\subsection{Localization}\label{subseclocaliz}

This section and the next one carry out Step 1 of the 
outline in Section \ref{subsecexoutline}. Specifically this 
Section shows that one can reduce 
from proving uniqueness in ${\cal W}$ to proving uniqueness 
in each of a family of bounded domains $\big\{{\cal W}_n\big\}_{n
\geq 0}$. 

By Theorem \ref{th:equiv}, uniqueness in distribution for 
$(\ref{eq:SDER})$ is equivalent to uniqueness for the 
constrained martingale problem for 
$(\mathbb{A},{\cal W},\mathbb{B},\Xi )$. Let 
\begin{equation}\begin{array}{c}
{\cal W}_0:={\cal W}\cap B_4,\quad\partial_j{\cal W}_0:=\partial_
j{\cal W}\cap B_4,\,j=1,\cdots ,m,\quad\partial_{m+1}{\cal W}_0:=
\partial B_4\cap {\cal W},\non\\
\hfill \\
g^{m+1}(x):=n^{m+1}(x)=-\frac x{|x|},\mbox{\rm \ for }x\in\partial 
B_4,\non\end{array}
\label{eq:loc0dom}\end{equation}
and $G_0(x)$, $\Xi_0$ be defined by $(\ref{eq:refcone})$, 
$(\ref{eq:refcone0})$ and $(\ref{eq:clset})$ with ${\cal W}$ replaced 
by ${\cal W}_0$, $g^j$, $j=1,\cdots ,m$ as in Condition \ref{G} (i) and 
$g^{m+1}$ as above.  
For $n\geq 1$, let 
\[{\cal W}_n:={\cal W}\cap\big(\overline {B_{1+2(n-1)}}\big)^c\cap 
B_{6+2(n-1)},\quad\partial_j{\cal W}_n:=\partial_j{\cal W}\cap\big
(\overline {B_{1+2(n-1)}}\big)^c\cap B_{6+2(n-1)},
,\]
\[\partial_0{\cal W}_n:=\partial B_{1+2(n-1)}\cap {\cal W},\quad\partial_{
m+1}{\cal W}_n:=\partial B_{6+2(n-1)}\cap {\cal W},\]
\[g^0(x):=n^0(x)=\frac x{|x|},\mbox{\rm \ for }x\in\partial B_{1+
2(n-1)},\]
 
\[g^{m+1}(x):=n^{m+1}(x)=-\frac x{|x|},\mbox{\rm \ for }
x\in\partial B_{6+2(n-1)},\]
and $G_n(x)$, $\Xi_n$ be defined by $(\ref{eq:refcone})$ and $(\ref{eq:clset})$ 
with ${\cal W}$ replaced 
by ${\cal W}_n$, $g^j$, $j=1,\cdots ,m$ as in Condition \ref{G} (i) and 
$g^0$, $g^{m+1}$ as above.

\begin{proposition}\label{th:localiz}
Assume Conditions \ref{G} and \ref{auxfunc}. 
If uniqueness holds for natural solutions of the 
constrained martingale problem for $(\A,{\cal W}_n,\B,\Xi_n)$ for 
each $n\geq 0$, then uniqueness holds for natural solutions of 
the constrained martingale problem for $(\A,{\cal W},\B,\Xi )$. 
\end{proposition}

\begin{proof}
See the Appendix.
\end{proof}

\subsection{Uniqueness in a domain ${\cal W}\cap\big(\overline {B_{
\delta}}\big)^c\cap B_{\delta'}$}\label{subsecDI}
Completing Step 1 of the outline in Section 
\ref{subsecexoutline}, let us show that uniqueness holds 
in each domain of the form 
\begin{equation}{\cal W}_{\delta ,\delta'}:={\cal W}\cap\big(\overline {
B_{\delta}}\big)^c\cap B_{\delta'}\qquad 0<\delta <\delta',\label{eq:locdeldom}\end{equation}
Let 
\[\partial_j{\cal W}_{\delta ,\delta'}:=:=\partial_j{\cal W}\cap\big
(\overline {B_{\delta}}\big)^c\cap B_{\delta'},\quad
\partial_0{\cal W}_{\delta ,\delta'}:=\partial B_{\delta}\cap {\cal W}
,\quad\partial_{m+1}{\cal W}_{\delta ,\delta'}:=\partial B_{\delta'}
\cap {\cal W},\]
\begin{equation}g^0(x):=n^0(x)=\frac x{|x|},\mbox{\rm \ for }x\in
\partial B_{\delta},\quad g^{m+1}(x):=n^{m+1}(x)=-\frac x{|x|},\mbox{\rm \ for }
x\in\partial B_{\delta'},\label{eq:locdelref}\end{equation}
and $G_{\delta ,\delta'}(x)$, $\Xi_{\delta ,\delta'}$ be defined by $
(\ref{eq:refcone})$ 
and $(\ref{eq:clset})$ with ${\cal W}$ 
replaced by ${\cal W}_{\delta ,\delta'}$, $g^j$, $j=1,\cdots ,m$ as in Condition \ref{G} (i) and 
$g^0$, $g^{m+1}$ as above.

\begin{lemma}\label{th:DI}
For every $0<\delta <\delta'$, uniqueness holds for natural solutions of the constrained 
martingale problem for  
$(\mathbb{A},{\cal W}_{\delta ,\delta'},\mathbb{B},\Xi_{\delta ,\delta'}
)$. 

In particular uniqueness holds for natural solutions of the 
constrained martingale problem for $(\A,{\cal W}_n,\B,\Xi_n)$ for 
each $n\geq 1$. 
\end{lemma}

\begin{proof}
By the analog of Theorem \ref{th:equiv}, it is equivalent to prove 
uniqueness in distribution for $(\ref{eq:SDER})$ in ${\cal W}_{\delta 
,\delta'}$ 
with cone of directions of reflection $G_{\delta ,\delta'}$. Conditions 
\ref{G} (i), (ii) and (iii) imply that the assumptions of 
Case 2 of \cite{DI93} are verified (see Remark 3.1 of 
\cite{DI93}). Therefore, by Corollary 5.2 of \cite{DI93}, 
the solution to $(\ref{eq:SDER})$ is pathwise unique, and 
hence unique in distribution. 
\end{proof}

The next lemma carries out the first part of Step 2 
of the outline in Section \ref{subsecexoutline}, namely 
proves that the distribution of a natural solution of 
the constrained martingale problem for 
$(\A,{\cal W}_0,\B,\Xi_0)$ is uniquley determined till the first time 
the solution hits the origin. The second statement in 
Step 2 that, starting off 
the origin, with probability one the origin is never 
reached, 
holds if Condition \ref{auxfunc} (i) is verified, as 
in the example of Section \ref{sectionex}. In this case 
Step 4 has a simpler proof. However Step 4 can be carried out 
as well if the origin can be reached (see Section 
\ref{subseczerotime0}).

For a solution $X$ of the constrained martingale problem for 
$(\A,{\cal W}_0,\B,\Xi_0)$ and $0<\delta <4$, let $\vartheta^X$ and $
\tau^{X,\delta}$ be defined by 
$(\ref{eq:0hittime})$ and 
$(\ref{eq:deltahittime})$ respectively. When there is no risk of confusion 
the superscript $X$ will be omitted. 

\begin{lemma}\label{th:unq-up-to-0}
For every initial distribution $\nu$ with support in 
${\cal W}_0-\{0\}$, $X(\cdot\wedge\vartheta^X)$ has the same distribution for all 
natural solutions $X$ of the constrained martingale problem for 
$(\A,{\cal W}_0,\B,\Xi_0)$ with initial distribution $\nu$.  
\end{lemma}

\begin{proof}
For every natural solution $X$ of the constrained martingale problem for 
$(\A,{\cal W}_0,\B,\Xi_0)$, 
\[X(\cdot\wedge\vartheta^X)=\lim_{k\rightarrow\infty}X(\cdot\wedge
\tau^{X,{\cal W}_0\cap\big(\overline {B_{1/k}}\,\big)^c}),\]
where $\tau^{X,{\cal W}_0\cap\big(\overline {B_{1/k}}\,\big)^c}:=\inf
\{t\geq 0:\,X(t)\in\overline {B_{1/k}}\,\}$. 

However 
$X(\cdot\wedge\tau^{X,{\cal W}_0\cap\big(\overline {B_{1/k}}\,\big
)^c})$ is a solution of the stopped 
constrained martingale problem for 
$(\mathbb{A},{\cal W}_{1/(2k),4},\mathbb{B},\Xi_{1/(2k),4};{\cal W}_{
1/(2k),4}\cap\big(\overline {B_{1/k}}\,\big)^c)$, 
which is unique by Lemma \ref{th:DI} and Theorem 2.11 
of \cite{CK24spa}. 
\end{proof}

\subsection{Exit time from $B_{\delta}$ and occupation time of the 
origin}\label{subseczerotime0}

This section carries out Step 3 of the outline in Section 
\ref{subsecexoutline}. 
As mentioned in Section \ref{subsecDI}, in general a natural solution of 
the constrained martingale problem for 
$(\A,{\cal W}_0,\B,\Xi_0)$ can reach the 
origin with positive probability. However Step 4 can 
still be carried out, provided one shows that every solution 
spends zero time at the origin.

\begin{lemma}\label{th:zerotime0}
There exists $\bar{\delta }>0$, $\bar {c}>0$, depending only on the data of 
the problem, such that, for $\delta\leq\bar{\delta}$, for every 
natural solution $X$ of the 
constrained martingale problem for 
$(\mathbb{A},{\cal W}_0,\mathbb{B},\Xi_0)$ starting at $0$, 
\[\E[\tau^{X,\delta}]\leq\bar {c}\delta^2.\]

For every natural solution $X$ of the constrained martingale 
problem for $(\mathbb{A},{\cal W}_0,\mathbb{B},\Xi_0)$, 
\begin{equation}\int_0^{\infty}{\bf 1}_{\{0\}}(X(t))\,dt=0,\quad 
a.s..\label{eq:zerotime0}\end{equation}
$(\ref{eq:zerotime0})$ implies in particular that, if $X$ starts at $
0$, 
\[\lim_{\delta\rightarrow 0}\tau^{X,\delta}=0,\qquad 
a.s..\]
\end{lemma}

\begin{proof}
The first assertion follows from Condition \ref{G} 
(iv) as in Section \ref{subsecexleave0}. The proof is 
analogous to that of Lemma 4.2 of \cite{CK18} and Lemma 
6.4 of \cite{TW93}.

For the second assertion, the proof of Lemma 3.21 of \cite{CK24aap} 
carries over. The argument is similar to that of Lemma 
2.1 of \cite{TW93}. It relies on the fact that  
every natural solution $X$ of the constrained martingale 
problem for $(\mathbb{A},{\cal W}_0,\mathbb{B},\Xi_0)$ is a 
continuous semimartingale with $<X>(t)=\sigma\sigma^Tt$ and on 
Tanaka's formula for the local time. 

The third assertion follows by observing that 
\begin{eqnarray*}
&&0\leq\lim_{\delta\rightarrow 0}\big(\tau^{X,\delta}\wedge 1\big)=\lim_{
\delta\rightarrow 0}\int_0^{\tau^{X,\delta}\wedge 1}{\bf 1}_{\overline {
B_{\delta}(0)}}(X(t))\,dt\\
&&\leq\lim_{\delta\rightarrow 0}
\int_0^1{\bf 1}_{\overline {B_{\delta}(0)}}(X(t))\,dt=\int_0^1{\bf 1}_{
\{0\}}(X(t))\,dt=0.
\end{eqnarray*}
\end{proof}

\subsection{Uniqueness and hitting distributions}\label{subsechitdist}

The first goal of this section is to prove the claim of 
Step 1 of the outline in Section \ref{subsecexoutline} 
that, in order to prove uniqueness among natural solutions of the 
constrained martingale problem for $(\A,{\cal W}_0,\B,\Xi_0)$ one can 
reduce to proving uniqueness among strong Markov, 
natural solutions (Lemma \ref{th:stM}). Next Step 4 of 
the outline is carried out (Lemma \ref{th:unq-by-hit}).

\begin{lemma}\label{th:stM}
There exist strong Markov natural solutions of the 
constrained martingale problem for $(\A,{\cal W}_0,\B,\Xi_0)$. 

If uniqueness holds among strong Markov solutions, then 
uniqueness holds for natural solutions of the constrained 
martingale problem for $(\A,{\cal W}_0,\B,\Xi_0)$. 
\end{lemma}

\begin{proof}
Existence of solutions of the controlled martingale 
problem for $(\A,{\cal W}_0,\B,\Xi_0)$ can be proved as in the first 
step of the proof of Theorem 3.23 of 
\cite{CK24aap}, hence Conditions 3.5 a) and b) of \cite{CK19}
are satisfied. 
By Lemma \ref{th:lambda0pos} and Lemma 3.3 of 
\cite{CK19}, Condition 3.5 c) of \cite{CK19} is also satisfied. 
Then the assertion follows from Corollaries 4.12 and 4.13 
of \cite{CK19}. 
\end{proof}

\begin{lemma}\label{th:unq-by-hit}
Suppose that, for each $\delta >0$ sufficiently small, the hitting 
distribution 
\[{\cal L}(X(\tau^{X,\delta}))(B)=\mathbb P\{X(\tau^{X,\delta})\in 
B\},\quad B\in {\cal B}\big(\partial B_{\delta}\big),\]
is the same for all strong Markov, natural solutions $X$ 
of the constrained martingale problem for 
$(\mathbb{A},{\cal W}_0,\mathbb{B},\Xi_0)$ starting at $0$. 

Then uniqueness holds among strong Markov natural solutions of the 
constrained martingale problem for $(\mathbb{A},{\cal W}_0,\mathbb{
B},\Xi_0)$. 
\end{lemma}

\begin{proof}
The proof of Lemma 3.27 of \cite{CK24aap} 
carries over. The argument is more complex but similar to that in 
Section \ref{subsecexhitdist}. Instead of 
$(\ref{eq:exappres})$, we now 
have, for any strong Markov natural solution of the 
constrained martingale problem for 
$(\mathbb{A},{\cal W}_0,\mathbb{B},\Xi_0)$, $X$,
\begin{eqnarray}
&&\mathbb E\left[\int_0^{\infty}e^{-\eta t}f(X(t))dt\right]\\
&&=\mathbb E\left[\int_0^{\vartheta}e^{-\eta t}f(X(t))dt\right]+\lim_{
n\rightarrow\infty}\mathbb E\left[\sum_{l=1}^{\infty}\prod_{m=0}^{
l-1}e^{-\eta (\vartheta_m^n-\tau_m^n)}\int_{\tau_l^n}^{\vartheta_
l^n}e^{-\eta (t-\tau_l^n)}f(X(t))dt\right],\non\label{eq:appres}\end{eqnarray}
where, for each $n\geq 0$,  
\begin{eqnarray*}
\vartheta_0^n:=\vartheta :=\inf\{t\geq 0:X(t)=0\},\qquad\qquad\qquad
\qquad\qquad\\
\tau_l^n:=\inf\{t\geq\vartheta_{l-1}^n:|X(t)|=\delta\,2^{-2n}\},\quad
\vartheta_l^n:=\inf\{t\geq\tau_l^n:X(t)=0\},\qquad l\geq 1,\end{eqnarray*}
and with the convention that $e^{-\infty}=0$. The limit holds by 
$(\ref{eq:zerotime0})$. 
\end{proof}

\subsection{The ergodic argument}\label{subsecergodic}

This section and Sections \ref{subsecLyap} and 
\ref{subsecscale} below carry out Step 5 of the outline in 
Section \ref{subsecexoutline}. By Section 
\ref{subsechitdist} we are reduced to proving that, for 
all $\delta >0$ sufficiently small, for any strong Markov, natural solution $
X$ 
of the constrained martingale problem for 
$(\mathbb{A},{\cal W}_0,\mathbb{B},$ $\Xi_0)$ starting at the origin, 
the hitting distribution on $\partial B_{\delta}$ is always the same. 

Let $X$ be a strong Markov, natural solution of the constrained martingale problem 
for $(\mathbb{A},{\cal W}_0,\mathbb{B},\Xi_0)$ starting at the origin 
and define, for $0<\delta <4$, 
\begin{equation}\tau^l:=\inf\{t\geq 0:|X(t)|=2^{-2l}\delta \},\quad 
l\in\Z_{+}.\label{eq:hittimes}\end{equation}
Note that, since $X$ starts at the origin, $\tau^l<\tau^{l-1}$. 
As discussed in Remark \ref{re:exkilled}, in the general 
case a solution can hit the origin 
between $\tau^l$ and $\tau^{l-1}$ with positive probability. Then, 
since the distribution of a solution starting on $\partial B_{2^{
-2l}\delta}$ 
is uniquely determined only until it hits the origin, 
in order to have transition kernels that are the 
same for all solutions, we cannot 
consider the transition kernels $(\ref{eq:extransker})$, but 
we have to consider the subprobability transition kernels 
\begin{equation}\begin{array}{c}
Q_{n-k}(x,B):=\P\big(\tau^{n-k-1}<\vartheta^{n-k},\,X(\tau^{n-k-1}
)\in B\big|X(\tau^{n-k})=x\big),\\
\hfill \\
B\in {\cal B}\big(\partial B_{2^{-2(n-k-1)}\delta}\big),\;x\in\partial 
B_{2^{-2(n-k)}\delta},\;0\leq k\leq n-1,\;n\in\N,\end{array}
\label{eq:transker}\end{equation}
where 
\begin{equation}\vartheta^{n-k}:=\inf\{t\geq\tau^{n-k}:\,X(t)=0\}
.\label{eq:0hittimes}\end{equation}
However, Lemma \ref{th:transker} below shows that one 
can generalize $(\ref{eq:exbasid})$ and still employ an ergodic type 
argument as in Section \ref{subsecexergodic}. 

\begin{lemma}\label{th:transker}

For $n\in\N$, $0\leq k\leq n-1$, let $Q_{n-k}$ be the integral operator 
corresponding to $(\ref{eq:transker})$, 
where $X$ is a strong Markov, natural solution of the constrained martingale problem 
for $(\mathbb{A},{\cal W}_0,\mathbb{B},\Xi_0)$ starting at the origin; 
$Q_{n-k}$ is uniquely defined by 
Lemma \ref{th:unq-up-to-0}. 

Then, for every strong Markov, natural solution $X$ 
of the constrained martingale problem for 
$(\mathbb{A},{\cal W}_0,\mathbb{B},\Xi_0)$ 
starting at $0$, 
\begin{equation}\E[f(X(\tau^{X,\delta}))]=\frac {\int Q_n\cdots Q_
1f(x)\nu_n(dx)}{\int Q_n\cdots Q_11(x)\nu_n(dx)},\quad\forall f\in 
{\cal C}\big(\partial B_{\delta}\big),\quad\forall n\geq 1,\label{eq:basid}\end{equation}
where $\nu_n$ is defined by 
\begin{equation}\nu_n(C):=\mathbb P\{X(\tau^n)\in B\},\quad B\in 
{\cal B}\big(\partial B_{2^{-2n}\delta}\big),\quad n\geq 0,\label{eq:hitdist}\end{equation}
and $\tau^n$ is as in $(\ref{eq:hittimes})$. 
\end{lemma}

\begin{proof}
The proof is the same as for Lemma 3.27 of 
\cite{CK24aap}. It relies on the strong Markov property. 
\end{proof}

Note that $(\ref{eq:basid})$ reduces to $(\ref{eq:exbasid})$ 
if 
\[Q_n\cdots Q_11(x)=\P\big(\tau^{\delta}<\vartheta^n\big|X(\tau^n
)=x\big)=1,\quad\forall x\in\partial B_{2^{-2n}\delta},\;\forall 
n\geq 1,\]
as in the case of the example of Section \ref{sectionex}. 

The ergodic theorem that is used to prove both Theorem 
\ref{th:exuniq} and Theorem \ref{th:uniq} 
is the following {\em reverse ergodic theorem for inhomogeneous }
{\em killed Markov chains}, which is proved in Section 2 of 
\cite{CK24aap}. 

\begin{theorem}\label{th:ergodic}
Let $E_0,E_1,E_2,\ldots$ be a sequence of 
compact metric spaces and $Q_1,Q_2,\ldots$ be a sequence of 
sub-probability transition kernels, with $Q_l$ governing 
transitions from $E_l$ to $E_{l-1}$. For $x,\tilde {x}\in E_l$, let 
\begin{equation}f_{l,\tilde {x}}(x,y):=\frac {dQ_l(x,\cdot )}{d\big
(Q_l(x,\cdot )+Q_l(\tilde {x},\cdot )\big)}(y)\label{eq:reldens}\end{equation}
and
\begin{equation}\epsilon_l(x,\tilde {x}):=\int\left(f_{l,\tilde {
x}}(x,y)\wedge f_{l,x}(\tilde x,y)\right)(Q_l(x,dy)+Q_l(\tilde {x}
,dy)).\label{eq:tvdensbnd}\end{equation}

Assume $\sup_xQ_l(x,E_{l-1})>0$, for all $l$, and there exist 
$c_0>0$ and $\epsilon_0>0$ such that 

\begin{itemize}
\item[(i)] 
\[\inf_{x\in E_n}Q_n\cdots Q_11(x)\geq c_0\sup_{x\in E_n}Q_n\cdots 
Q_11(x),\quad\forall n\geq 1,\]
\item[(ii)] 
\[\inf_n\inf_{x,\tilde {x}\in E_n}\epsilon_n(x,\tilde {x})\geq\epsilon_
0.\]
\end{itemize}
Then 
\[\sup_{x\in E_n}Q_n\cdots Q_11(x)>0,\quad\forall k,\]
and, for any sequence $\{\nu_n\}$, $\nu_n$ a probability measure on 
$E_n$, 
\[\lim_{n\rightarrow\infty}\frac {\int Q_nQ_{n-1}\cdots Q_1f(x)\,
\nu_n(dx)}{\int Q_nQ_{n-1}\cdots Q_11(x)\,\nu_n(dx)}=C(f),\qquad\forall 
f\in {\cal C}(E_0),\]
where $C(f)$ is independent of the sequence $\{\nu_n\}$. 
\end{theorem}\

Sections \ref{subsecLyap} and \ref{subsecscale} below show that, 
under Conditions \ref{G} and \ref{auxfunc}, 
assumptions (i) and (ii), respectively, are verified by the 
transition kernels $(\ref{eq:transker})$ for $\delta$ sufficiently small, 
so that Theorem \ref{th:ergodic} can be applied. Then the fact that, 
for any strong Markov, natural solution $X$ 
of the constrained martingale problem for 
$(\mathbb{A},{\cal W}_0,\mathbb{B},\Xi_0)$ starting at the origin, 
the hitting distribution on $\partial B_{\delta}$ is always the same will 
follow from Lemma \ref{th:transker}. 

Finally, the proof of Theorem \ref{th:uniq} is 
summarized at the end of Section \ref{subsecscale}.

\subsection{Uniform bound on hitting times}\label{subsecLyap}

For the transition kernels $(\ref{eq:transker})$, assumption (i) of Theorem 
\ref{th:ergodic} translates into the following uniform 
lower bound. Let $X$ be a 
strong Markov, natural solution 
of the constrained martingale problem for 
$(\mathbb{A},{\cal W}_0,\mathbb{B},\Xi_0)$ starting at the origin 
and let $\tau^{\delta},$ $\tau^n$ and $\vartheta^n$ be defined by $
(\ref{eq:deltahittime})$, 
$(\ref{eq:hittimes})$ and $(\ref{eq:0hittimes})$ (with $n-k$ 
replaced by $n$) respectively.
Supposing 
\[\P\big(\tau^{X,\delta}<\vartheta^n\big|X(\tau^n)=x\big)>0,\quad
\forall x\in\overline {{\cal W}}\cap\partial B_{2^{-2n}\delta},\;
\forall n\geq 1,\]
there exists $c_0>0$ such that 
\begin{equation}\inf_{n\geq 1}\;\inf_{x,\tilde {x}\in\overline {{\cal W}}
\cap\partial B_{2^{-2n}\delta}}\frac {\P\big(\tau^{X,\delta}<\vartheta^
n\big|X(\tau^n)=x\big)}{\P\big(\tau^{X,\delta}<\vartheta^n\big|X(
\tau^n)=\tilde {x}\big)}\geq c_0.\label{eq:hit-tbound}\end{equation}
Recall that, by Lemma \ref{th:unq-up-to-0}, the 
numerator in $(\ref{eq:hit-tbound})$ 
is the same for all strong Markov, natural solutions of 
the constrained martingale problem for 
$(\mathbb{A},{\cal W}_0,\mathbb{B},\Xi_0)$ 
starting at $0$ and is the same as 
\[\P\big(\tau^{X^x,\delta}<\vartheta^{X^x}\big),\]
where $X^x$ is a solution starting at $x\in\overline {{\cal W}}$, $
|x|=2^{-2n}\delta$, and 
analogously for the denominator.

\begin{lemma}\label{th:hitprob}
\begin{itemize}
\item[(i)]If Condition \ref{auxfunc} (i) is verified, 
\begin{equation}\P\big(\tau^{X^x,\delta}<\vartheta^{X^x}\big)=1,\qquad\forall x
\in\overline {{\cal W}},\;0<|x|<\delta\leq\delta_{{\cal W}}.\label{eq:0nothit}\end{equation}
\item[(ii)]If Condition \ref{auxfunc} (ii) is verified, 
\begin{equation}\P\big(\tau^{X^x,\delta}<\vartheta^{X^x}\big)>0,\qquad\forall x
\in\overline {{\cal W}},\;0<|x|<\delta\leq\delta_{{\cal W}}.\label{eq:0poshit}\end{equation}
Moreover there exists $c_0>0$ such that, for all 
$0<\delta\leq\delta_{{\cal W}}$,  
\begin{equation}\inf_{x,\tilde {x}\in\overline {{\cal W}}:\,0<|x|
=|\tilde {x}|<\delta}\frac {\P\big(\tau^{X^x,\delta}<\vartheta^{X^
x}\big)}{\P\big(\tau^{X^{\tilde {x}},\delta}<\vartheta^{X^x}\big)}
\geq c_0.\label{eq:0hitbnd}\end{equation}
\end{itemize}
Of course, if $(\ref{eq:0nothit})$ holds $(\ref{eq:0poshit})$ 
and $(\ref{eq:0hitbnd})$ hold too. 
\end{lemma}

\begin{proof}
The proof of Lemma 3.29 of \cite{CK24aap} carries over. 

If Condition \ref{auxfunc} (i) is verified, the assertion 
follows as in Section \ref{subsecex0nothit} 
by applying Ito's formula to the function $V$ between times 
$0$ and $\tau^{X,\epsilon}\wedge\tau^{X,\delta}\wedge t$, $0<\epsilon 
<|x|$, and taking limits as $t\rightarrow\infty$ first, and 
$\epsilon\rightarrow 0$ next. 

If Condition \ref{auxfunc} (ii) is verified, the first assertion 
follows by applying Ito's formula to the 
functions $V_{+}$ an $V_{-}$ between the same times as above and 
taking the same limits, but using $(\ref{eq:V+0})$. 
The second assertion then follows by $(\ref{eq:V+unif})$. 
\end{proof}

\subsection{Uniform bound on hitting distributions of 
$\partial B_{\delta}$ and proof of Theorem \ref{th:uniq}}\label{subsecscale}

Both for the example of Section \ref{sectionex} and in the 
general case, in order to be able to apply Theorem 
\ref{th:ergodic} we still need to verify the lower bound 
$(ii)$. This is done in the following two lemmas. 
Finally, at the end of this section, the proof of Theorem \ref{th:uniq} is 
summarized.

\begin{lemma}\label{th:scaling}
For any sequence $\{x^n\}\subseteq\overline {{\cal W}}$ such that $
\{2^{2n}x^n\}$ converges to 
some $\bar {x}\in\overline {{\cal W}}-\{0\}$, let $X^{x^n}$ be a natural solution of the constrained 
martingale problem for $(\mathbb{A},{\cal W},\mathbb{B},\Xi )$ starting at $
x^n$.

Then the sequence of processes $\bigg\{2^{2n}\,X^{x^n}(2^{-4n}\cdot 
)\bigg\}$ is 
relatively compact and any of its limit points, $\bar {X}^{\bar {
x}}$, 
is a solution of $(\ref{eq:SDER})$ with $b=0$ and 
$\bar {X}^{\bar {x}}(0)=\bar {x}$. 

In particular, for every open set ${\cal O}$ such that 
${\cal O}\cap {\cal W}\cap\partial B_{2\delta}(0)\neq\emptyset$, there exists a constant $
\eta_0>0$ such 
that, for $|x^n|=2^{-2n}\delta$, $\{2^{2n}x^n\}$ converging to $\bar {
x}$, and ${\cal O}^n:=\{x:\,2^{2n}x\in {\cal O}\}$, 
\begin{equation}\liminf_{n\rightarrow\infty}\P (\sigma^n<\vartheta^{
X^{x^n}},\,\,\,X^{x^n}(\sigma^n)\in {\cal O}^n)\geq\eta_0,\label{eq:pos-hit}\end{equation}
where 
\[\sigma^n:=\inf\{t\geq 0:|X^{x^n}(t)|=2^{-2n+1}\delta \}.\]
\end{lemma}

\begin{proof}
The first two assertions can be proved by the same 
time change and compactness arguments as in Lemma 4.5 
and Theorem 4.1 of \cite{CK18}. 

Let us prove the last assertion. As in Lemma 3.30 of \cite{CK24aap}, let 
$\psi :[0,t_0]\rightarrow\R^d$ be a continuous function such that 
\begin{eqnarray*}
\psi (0)=\bar {x},\qquad |\psi (s)|=|\bar {x}|=\delta ,\;\psi (s)
\in {\cal W},\mbox{\rm \ for }0<s\leq\frac {t_0}2,\qquad 2\psi (\frac {
t_0}2)\in {\cal O}\cap {\cal W},\\
\frac {\psi (s)}{|\psi (s)|}=\frac {\psi (\frac {t_0}2)}{|\psi (\frac {
t_0}2)|},\quad |\psi (s)|=\delta\frac 2{t_0}(t_0-s)+(2\delta +\frac {
\epsilon}2)\frac 2{t_0}(s-\frac {t_0}2),\mbox{\rm \ for }\frac {t_
0}2<s\leq t_0\end{eqnarray*}
and 
\[\epsilon <\delta ,\quad B_{\epsilon}(2\psi (\frac {t_0}2))\subseteq 
{\cal O}\cap {\cal W}.\]
Then 
\begin{eqnarray*}
&&\liminf_{n\rightarrow\infty}\P\big(\sigma^n<\vartheta^{X^{x^n}}
,\,\,\,X^{x^n}(\sigma^n)\in {\cal O}^n\big)\geq\liminf_{n\rightarrow
\infty}\P\big(\sup_{0\leq s\leq t_0}\big|X^{x^n}(s)-\psi (s)\big|
<\frac {\epsilon}2\big)\\
&&\geq\P\big(\sup_{0\leq s\leq t_0}\big|\bar {X}^{\bar {x}}(s)-\psi 
(s)\big|<\frac {\epsilon}2\big)=P^{\bar {X}^{\bar {x}}}\big\{x\in 
{\cal C}_{[0,\infty )}(\overline {{\cal W}}):\,\big|x(s)-\psi (s)\big
|<\frac {\epsilon}2\big\},\end{eqnarray*}
where $P^{\bar {X}^{\bar {x}}}$ is the law of $\bar {X}^{\bar {x}}$. 
Since $\delta\leq |\psi (s)|$ for all $0\leq s\leq t_0$, the right hand side equals 
\[P^{\bar {X}^{\bar {x}}(\cdot\wedge\vartheta^{\bar {X}^{\bar {x}}}
)}\big\{x\in {\cal C}_{[0,\infty )}(\overline {{\cal W}}):\,\big|
x(s)-\psi (s)\big|<\frac {\epsilon}2\big\},\]
where $P^{\bar {X}^{\bar {x}}(\cdot\wedge\vartheta^{\bar {X}^{\bar {
x}}})}$ is the law of $\bar {X}^{\bar {x}}(\cdot\wedge\vartheta^{
\bar {X}^{\bar {x}}})$ and is 
uniquely determined by Lemma \ref{th:unq-up-to-0}. In Lemma 3.30 of 
\cite{CK24aap} the assertion then follows from Theorem 3.1 of 
\cite{KW91}. Here that theorem cannot be applied 
immediately, because the cone is not smooth. 
However a careful inspection of Theorem 3.1 of 
\cite{KW91} shows that it 
extends to a cone with piecewise smooth boundary, provided that 
\[\E\bigg[\int_0^{\vartheta^{\bar {X}^{\bar {x}}}}{\bf 1}_{\partial 
{\cal W}-\{0\}}(\bar {X}^{\bar {x}}(t))\,dt\bigg]=0.\]
However 
\begin{eqnarray*}
\E\bigg[\int_0^{\vartheta^{\bar {X}^{\bar {x}}}}{\bf 1}_{\partial 
{\cal W}-\{0\}}(\bar {X}^{\bar {x}}(t))\,dt\bigg]&=&\lim_{n\rightarrow
\infty}\E\bigg[\int_0^{\tau^{\bar {X}^{\bar {x}},\,\delta /n}}{\bf 1}_{
\partial {\cal W}}(\bar {X}^{\bar {x}}(t))\,dt\bigg]\\
&=&\lim_{n\rightarrow\infty}\E\bigg[\int_0^{\tau^{\tilde {X}^{n,\bar {
x}},\,\delta /n}}{\bf 1}_{\partial {\cal W}}(\tilde {X}^{n,\bar {
x}})(t))\,dt\bigg],\end{eqnarray*}
where $\tilde {X}^{n,\bar {x}}$ is the natural solution of the constrained martingale 
problem for $(\A,{\cal W}_{\delta /(2n),4},\B,$ $\Xi_{\delta /(2n)
,4})$, with $b=0$ and ${\cal W}_{\delta /(2n),4}$, $\Xi_{\delta /
(2n),4}$ defined by 
$(\ref{eq:locdeldom})$, $(\ref{eq:locdelref})$, 
$(\ref{eq:clset})$, 
and the last term on the right hand is zero by Proposition 2.12 of 
\cite{KR17} (note that every solution of a constrained 
martingale problem is a solution of the corresponding 
submartingale problem).
\end{proof}

\begin{lemma}\label{th:coupling}
The transition kernels $(\ref{eq:transker})$ verify 
assumption (ii) of Theorem \ref{th:ergodic}. 
\end{lemma}

\begin{proof}
The proof is the same as for Lemma 3.31
of \cite{CK24aap}. 
Let $x,\tilde {x}\in\partial B_{2^{-2n}\delta}$. Lemma 5.3 of \cite{CK18} and 
$(\ref{eq:pos-hit})$ allow to 
construct two strong Markov, natural solutions of the 
constrained martingale problem for 
$(\mathbb{A},{\cal W}_0,\mathbb{B},\Xi_0)$, starting at $x$ and 
$\tilde {x}$ respectively, that couple before reaching $\partial 
B_{2^{-2(n-1)}\delta}$.
This yields that 
\[\|Q_n(x,\cdot )-Q_n(\tilde {x},\cdot )\|_{TV}\leq Q_n(x,\partial 
B_{2^{-2(n-1)}\delta})\vee Q_n(\tilde {x},\partial B_{2^{-2(n-1)}
\delta})-\epsilon_0,\]
for some positive constant $\epsilon_0$. On the other hand 
\[Q_n(x,\partial B_{2^{-2(n-1)}\delta})\vee Q_n(\tilde {x},\partial 
B_{2^{-2(n-1)}\delta})\leq\|Q_n(x,\cdot )-Q_n(\tilde {x},\cdot )\|_{
TV}+\epsilon_n(x,\tilde {x}).\]
\end{proof}
\vskip.2in

\noindent
{\bf Proof of Theorem \ref{th:uniq}}
By Theorem \ref{th:localiz} and Lemma \ref{th:DI} it is 
sufficient to prove uniqueness for natural solutions of 
the constrained martingale problem for 
$(\mathbb{A},{\cal W}_0,\mathbb{B},\Xi_0)$. To this end, by Lemma \ref{th:stM} 
and Lemma \ref{th:unq-by-hit}, it is enough to show 
that the hitting distribution on $\partial B_{\delta}$, ${\cal L}
(X(\tau^{\delta}))$, is the same 
for all strong Markov, natural solutions of the 
constrained martingale probem for 
$(\mathbb{A},{\cal W}_0,\mathbb{B},\Xi_0)$ starting at the origin. 
By Lemma \ref{th:transker}, this will be proved if we 
can apply Theorem \ref{th:ergodic} to the transition 
kernels $(\ref{eq:transker})$. 

In order to be able to apply Theorem \ref{th:ergodic}, 
one needs to verify the two assumptions (i) and (ii). 
Assumption (i) is verified by Lemma \ref{th:hitprob}, 
while assumption (ii) is verified by Lemma 
\ref{th:coupling}, 
\hfill $\Box$
\medskip

\appendix\section{Appendix}
\renewcommand {\thesection}{A}
\setcounter{equation}{0} 
\vskip.1in

\noindent {\bf Proof of Theorem \ref{th:localiz}}
\vskip.1in
The proof consists essentially in verifying the 
assumptions of Theorem 2.13 (ii) of \cite{CK24spa}. 

Let 
\[U_0=:\overline {{\cal W}}\cap B_3,\qquad U_n:=\overline {{\cal W}}
\cap\big(\overline {B_{2+2(n-1)}}\big)^c\cap B_{5+2(n-1)},\quad n
\geq 1.\]
A natural solution of the stopped constrained 
martingale problem for $(\A,{\cal W},\B,$$\Xi ;U_n)$ is a solution 
that can be obtained by a solution of 
the stopped controlled martingale problem for 
$(\A,{\cal W},\B,$$\Xi ;U_n)$ by a time change analogous to that in 
\ref{def:nat}, that is 
the process first is stopped and then time changed (see 
Definitions 2.1, 2.4 and 2.5 of \cite{CK24spa}). The time 
change does not always produce a solution of the 
stopped constrained martingale problem, in fact this is 
one of the requirements of Theorem 2.13 (ii) of 
\cite{CK24spa}. Theorem 2.13 (ii) of \cite{CK24spa} requires the following: 
for every solution $(Y,\lambda_0,\Lambda_1)$ of the controlled martingale 
problem for $(\A,{\cal W},\B,\Xi )$, $\lambda_0(t)>0$ for all $t>
0$; 
for every solution $(Y^{U_n},\lambda_0^{U_n},\Lambda_1^{U_n})$ of the stopped controlled 
martingale problem for $(\A,{\cal W},\B,\Xi ;U_n)$, $Y^{U_n}\circ\big
(\lambda_0^{U_n}\big)^{-1}$
is a natural solution of the stopped constrained 
martingale problem for $(\A,{\cal W},\B,$$\Xi ;U_n)$; 
uniqueness holds for natural solutions of the stopped constrained 
martingale problem for $(\A,{\cal W},\B,\Xi ;U_n)$. 
Under these conditions, uniqueness 
holds for natural solutions of the constrained 
martingale problem for $(\A,{\cal W},\B,\Xi )$. 

The first assumption is verified 
by Lemma \ref{th:lambda0pos}. 

As far as the second assumption is concerned, 
let $(Y^{U_n},\lambda_0^{U_n},\Lambda_1^{U_n})$ be a solution of the stopped controlled martingale 
problem for $(\A,{\cal W},\B,\Xi ;U_n)$ with initial distribution 
supported in $\overline {{\cal W}_n}$. Then it can be immediately verified 
that $(Y^{U_n},\lambda_0^{U_n},\Lambda_1^{U_n})$ is a solution of the stopped 
controlled martingale problem for $(\A,{\cal W}_n,\B,\Xi_n;U_n)$. 
By Lemma \ref{th:lambda0pos} and Lemmas 3.3, 
3.4 and Corollary 3.9 b) of \cite{CK19}, $Y^{U_n}\circ\big(\lambda_
0^{U_n}\big)^{-1}$ is 
a natural solution of the stopped constrained martingale 
problem for $(\A,{\cal W}_n,\B,\Xi_n;U_n)$ and hence 
a natural solution of the stopped constrained martingale 
problem for $(\A,{\cal W},\B,\Xi ;U_n)$. 
If $(Y^{U_n},\lambda_0^{U_n},\Lambda_1^{U_n})$ is a solution of the stopped controlled martingale 
problem for $(\A,{\cal W},\B,\Xi ;U_n)$ with initial distribution not 
supported in $\overline {{\cal W}_n}$, let 
\[\tilde {Y}^{U_n}(t):=Y^{U_n}(t)\mathbf{1}_{\overline {{\cal W}_
n}}(Y^{U_n}(0))+y_n^0\mathbf{1}_{\overline {{\cal W}}\,-\,\overline {
{\cal W}_n}}(Y^{U_n}(0)),\quad t\geq 0,\]
\[\big(\tilde{\lambda}_0^{U_n},\tilde{\Lambda}_1^{U_n}\big):=\big
(\tilde{\lambda}_0^{U_n},\tilde{\Lambda}_1^{U_n}\big)\mathbf{1}_{\overline {
{\cal W}_n}}(Y_n(0)),\]
where $y_n^0$ is some fixed point in $\overline {{\cal W}_n}$. Then 
$(\tilde {Y}^{U_n},\tilde{\lambda}_0^{U_n},\tilde{\Lambda}_1^{U_n}
)$ is a solution of the stopped 
controlled martingale problem for $(\A,{\cal W}_n,\B,\Xi_n;U_n)$ with 
initial distribution supported in $\overline {{\cal W}_n}$  
and $\tilde {Y}^{U_n}\circ\big(\tilde{\lambda}_0^{U_n}\big)^{-1}$ is a natural solution of the 
stopped constrained martingale problem for $(\A,{\cal W},\B,\Xi ;
U_n)$. 
Since $Y^{U_n}\circ\big(\lambda_0^{U_n}\big)^{-1}(t)=Y^{U_n}(0)$ for all $
t\geq 0$ if 
$Y^{U_n}(0)\in\overline {{\cal W}}\,-\,\overline {{\cal W}_n}$, and $
Y^{U_n}\circ\big(\lambda_0^{U_n}\big)^{-1}(t)=\tilde {Y}^{U_n}\circ\big
(\tilde{\lambda}_0^{U_n}\big)^{-1}(t)$ 
for all $t\geq 0$ if $Y^{U_n}(0)\in\overline {{\cal W}_n}$, $Y^{U_
n}\circ\big(\lambda_0^{U_n}\big)^{-1}$ is indeed a 
a natural solution of the stopped constrained 
martingale problem for $(\A,{\cal W},\B,\Xi ;U_n)$. 

Finally, for every natural solution $X^{U_n}$ of the stopped constrained martingale 
problem for $(\A,{\cal W},\B,\Xi ;U_n)$, setting 
\[\tilde {X}^{U_n}(t)=X^{U_n}(t)\mathbf{1}_{\overline {{\cal W}_n}}
(X^{U_n}(0))+x_n^0\mathbf{1}_{\overline {{\cal W}}\,-\,\overline {
{\cal W}_n}}(X^{U_n}(0)),\quad t\geq 0,\]
where $x_n^0$ is a fixed point in $\overline {{\cal W}_n}$, the distribution of $
X^{U_n}$ 
is uniquely determined by its initial distribution and by 
the distribution of $\tilde {X}^{U_n}$ and $\tilde {X}^{U_n}$ is a natural solution 
of the stopped constrained martingale problem for 
$(\A,{\cal W}_n,\B,\Xi_n;U_n)$. On the other hand, 
Lemma \ref{th:lambda0pos} and Proposition 2.9 of 
\cite{CK24spa} yield that $(\A,{\cal W}_n,\B,\Xi_n)$
and $U_n$ verify Condition 2.7 of \cite{CK24spa}. Therefore, 
by Corollary 2.12 of \cite{CK24spa}, uniqueness holds for 
natural solutions of the stopped constrained martingale 
problem for $(\A,{\cal W}_n,\B,\Xi_n;U_n)$. 
\hfill $\Box$
\vskip.5in



\end{document}